\documentclass{article}
\usepackage{amssymb}
\usepackage{amsmath}
\usepackage{enumerate}
\usepackage{graphicx}
\usepackage{amscd}
\def\div{\mathfrak{Div}}

\title{\huge The space of complete embedded maximal surfaces with isolated singularities in the 3-dimensional Lorentz-Minkowski space $\l^3$}
\author{\Large Isabel Fern\'{a}ndez    \thanks{Research partially
supported by
MEC-FEDER grant number MTM2004-00160. \newline
$\quad^{\dagger}$  Research partially supported by Consejería de Educaci\'{o}n y Ciencia de la Junta de Andaluc\'{\i}a and the European Union. \newline
}
\and \Large  Francisco J. L\'{o}pez $ ^{{\ast} {\dagger}}$ \and \Large  Rabah Souam$^{\dagger}$}

\usepackage{latexsym}
\usepackage{amssymb}
\usepackage{amsmath}
\usepackage{epsfig}
\usepackage[latin1]{inputenc}

\newcommand{\df}{ \stackrel{\rm def}{=}}

\def\h{\mathbb{H}}
\def\r{\mathbb{R}}
\def\n{\mathbb{N}}
\def\b{\mathbb{B}}
\def\c{\mathbb{C}}

\def\s{\mathbb{S}}

\def\l{\mathbb{L}}
\def\z{\mathbb{Z}}

\def\pg{\mathfrak{p}}
\def\Cg{\mathfrak{C}}
\def\Gg{\mathfrak{G}}

\def\sg{\mathfrak{s}}

\def\Jg{\mathfrak{J}}
\def\mg{\mathfrak{m}}
\def\vg{\mathfrak{v}}
\def\Mg{\mathfrak{M}}

\newenvironment{proof}{\trivlist
\item[\hskip\labelsep{\em Proof}\,:]}{\hfill{$\Box$}\endtrivlist}

\newtheorem{lemma}{Lemma}[section]
\newtheorem{remark}{Remark}[section]
\newtheorem{theorem}{Theorem}[section]
\newtheorem{proposition}{Proposition}[section]
\newtheorem{corollary}{Corollary}[section]

\newtheorem{definition}{Definition}[section]
\begin{document}
\maketitle

\begin{abstract}
We prove that a complete embedded maximal surface in $\l^3$ with a finite number of singularities is 
an entire maximal graph with conelike singularities over any spacelike plane, and so, it is asymptotic to a spacelike plane or a half catenoid.

We show that the moduli space $\Gg_n$ of entire maximal graphs over
$\{x_3=0\}$ in
$\l^3$ with $n+1 \geq 2$  singular points and vertical limit normal
vector at infinity is a $3n+4$-dimensional differentiable manifold. 
The convergence in $\Gg_n$ means the one of conformal structures and Weierstrass data, and it is equivalent to the uniform convergence of graphs on compact subsets of $\{x_3=0\}.$  Moreover, the position of the singular points in $\r^3$ and
the logarithmic growth at infinity can be used as global analytical coordinates with the same underlying topology.

We also introduce the moduli space ${\mathfrak{M}}
_n$ of {\em marked} graphs with $n+1$ singular points (a mark in a
graph is an ordering of its singularities), which is a $(n+1)$-sheeted covering of $\Gg_n.$ We prove that identifying  marked graphs  differing by translations,
rotations about a  vertical axis,
homotheties or symmetries about a horizontal  plane,  the corresponding
quotient space $\hat{\Mg}_n$ is an analytic manifold of dimension $3n-1.$ This manifold can be identified
with a spinorial bundle ${\cal S}_n$ associated  to the moduli
space of Weierstrass data of graphs in $\Gg_n.$
\end{abstract}

\section{Introduction} \label{sec:intro}

A maximal hypersurface in a Lorentzian manifold is a spacelike
hypersurface with zero mean curvature. Besides of their mathematical
interest these hypersurfaces and more generally those having constant
mean curvature have a significant importance  in classical Relativity.
More information on this aspect can be found for instance in
\cite{marsden-tipler}. When the ambient space is  the  Minkowski
space $\l^{n+1}$, one of the most  important
results  is the proof of a  Bernstein-type theorem  for
maximal hypersurfaces in $\l^{n+1}$.   Calabi
\cite{calabi}
$(n=2,3)$,  and Cheng and Yau \cite{cheng-yau} (for any $n$) have proved
that a complete  maximal hypersurface in $\l^{n+1}$ is necessarily a
spacelike hyperplane. It is therefore meaningless to consider global
problems
on maximal and {\em everywhere  regular}  hypersurfaces in $\l^{n+1}$.
Problems of interest should  deal with hypersurfaces having non empty
boundary or having a certain type of singularities. For instance,
Bartnik and Simon
\cite{Bar-Sim},  have obtained results on the existence and
regularity of (spacelike) solutions to the boundary value
problem for the mean curvature operator  in
$\l^{n+1}.$ Klyachin and Miklyukov \cite{kla-mik} have given results on
the existence of solutions, with a finite number of isolated
singularities, to the maximal hypersurface equation in $\l^{n+1}$ with
prescribed boundary conditions. 

More recently Klyachin \cite{kly} has studied the existence, uniqueness and asymptotic behavior of {\em entire} maximal graphs in $\l^{n+1}$ with prescribed flux vector at infinity and compact
singular set. A spacelike surface in $\l^3$ is said to be an entire graph if its orthogonal projection  over any spacelike plane is a homeomorphism. In particular, he proves that an entire maximal graph in $\l^3$ with a finite set of singularities is asymptotic to a half catenoid or a plane, and it is uniquely determined  by the flux vector at infinity and the values on the singular points. 

From a different point of view, Umehara and Yamada \cite{um-ya} have proved some results about the global behavior of maximal immersed
surfaces in $\l^3$ with analytical curves of singularities. They introduce the concept of complete {\em maxfaces} and prove an Osserman type inequality for this family of surfaces.
Complete embedded maximal surfaces with isolated singularities can be reflected analytically about singular points, leading to complete maxfaces in the sense of \cite{um-ya}. 

Maximal surfaces in $\l^3$ share some properties with minimal surfaces in the
Euclidean space $\r^3$. Both families  arise as solutions of variational
problems: local maxima (minima)
for the area functional in the Lorentzian (Euclidean) case. Like minimal
surfaces in $\r^3,$ the Gauss map is conformal and they 
admit a Weierstrass
representation in terms of meromorphic data (cf.\cite{kobayashi}).

For several reasons, isolated singularities of maximal surfaces in $\l^3$ are specially interesting. If the surface is embedded around the singularity, it is  locally a graph and the singular point is said to be of {\em conelike} type. Conelike singularities correspond to points where the Gauss curvature blows up, the Gauss map has no well defined limit and the surface is asymptotic to a half light cone (we refer to Ecker \cite{ecker}, Kobayashi \cite{kobayashi} and Klyachin and Miklyukov \cite{kly-mik} for a good setting).

\begin{figure}[htbp] 
\begin{center}
\includegraphics[width=14cm]{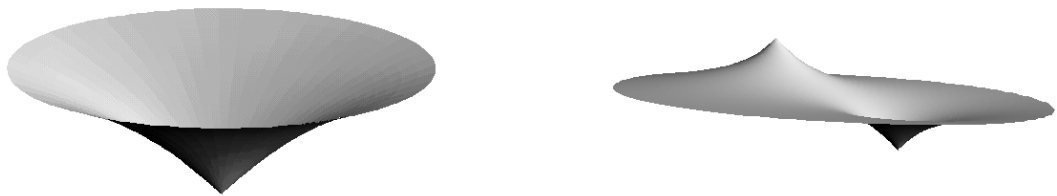}
\caption{The Lorentzian catenoid and a Riemann type example. \label{fig:cat-rie}}
\end{center}
\end{figure}

The simplest example of this phenomenon is the
{\it Lorentzian half catenoid} which is the only entire maximal graph over spacelike planes with an isolated singularity
( cf. \cite
{ecker}, and see \cite{kobayashi} for a
previous characterization). New examples of Riemann's type were  discovered
more recently  \cite{lopez-lopez-souam}. They form a one parameter family of  entire
maximal
graphs over spacelike planes asymptotic to a plane at infinity and having two isolated singularities. Moreover, they are characterized by the property of being foliated by complete circles, straight lines or singular points in parallel spacelike planes, besides the catenoid. In Section \ref{sec:global} we show explicit examples of entire maximal graphs with an arbitrary number of singular points. 

A maximal surface in $\l^3$ reflects about conelike points to its mirror, and its Weierstrass data extend by Schwarz reflection to the double. Following Osserman \cite{osserman}, this fact establishes an important connection between the theory of complete embedded maximal surfaces in $\l^3$ with a finite number of singularities  and classical algebraic geometry. Therefore, moduli problems arise in a natural way. 

This paper is  devoted to exploit this idea. We first observe
that a complete embedded maximal surface with a finite number of singularities is an entire maximal graph over any spacelike plane with conelike singularities, and so it is
asymptotic to a catenoid or a spacelike plane.  The corresponding moduli space
has  a structure of a finite dimensional real analytic  manifold and we
compute its dimension. The underlying topological structure correspond to the convergence of conformal structures and Weierstrass data, and it is equivalent to the uniform convergence of graphs on compact
subsets. 

Finite dimensional smoothness results for
the moduli spaces of other noncompact geometric objects have been
obtained in the past few years. Important contributions  are due to
Perez and Ros \cite{per-ros} in the case of non degenerate properly
embedded minimal surfaces with finite total curvature and fixed topology
in $\r^3,$ and to Kusner, Mazzeo and Pollack
\cite{kus-maz-pol} in the case of properly embedded non minimal
constant mean curvature surfaces in $\r^3$ with finite topology.  We
use different techniques to obtain our result. The approach of these
authors is analytic while ours takes advantage of the Weierstrass
representation and relies more on algebraic geometry tools (more precisely compact
Riemann surface theory).
\begin{figure}[htbp] 
\begin{center}
\includegraphics[width=14cm]{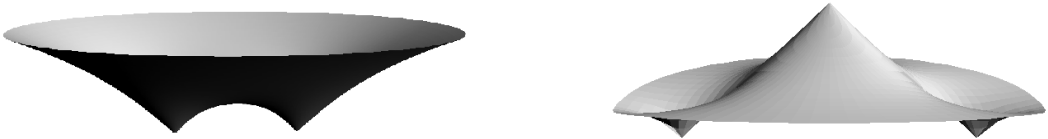}
\caption{Two graphs of catenoidal and  planar type with two and three singularities, respectively. \label{fig:2-3}}
\end{center}
\end{figure}

Our paper is organized as follows: 

Section
\ref{sec:prelim} is devoted to some preliminary results.
We review the local behavior of maximal surfaces around isolated singularities and the global geometry of
complete embedded maximal surfaces with a finite number of singular points,   stating these results
in terms of the Weierstrass data. In particular, we show that a complete embedded maximal surface with a finite number of singularities is an entire graph, has finite conformal type and its Weierstrass data extend meromorphically, in a controlled way, to its unique end. Moreover we characterize these surfaces in terms of meromorphic data on compact Riemann surfaces admitting a mirror involution, and use this result to construct examples with $n+1$ conical singularities for
any $n\geq 1.$ Finally, we  prove the following characterization of the Lorentzian catenoid:

\begin{quote}{\em The catenoid is the unique entire  maximal graph with
downward pointing conelike singularities and vertical flux at the singularities.}
\end{quote}

In Section \ref{sec:teich}, we define some natural bundles on the 
moduli space ${\cal T}_n$ of  marked unbounded planar circular domains with $n+1$
boundary components (the mark refers to an ordering of the boundary circles). We introduce a spinorial bundle ${\cal S}_n$ associated  to the moduli
space of Weierstrass data of surfaces in the space of graphs with $n+1$
conelike singularities (see Definiton \ref{def:cale}), and make the
fundamental analysis. In  Section \ref{sec:moduli} we prove
the main results  of the paper:

\begin{quote}{\em  The space $\Mg_n$ of marked entire maximal 
graphs over $\{x_3=0\}$ (the mark is an ordering of the set of singularities) with vertical limit normal vector at infinity and $n+1 \geq 2$  conelike singularities  in
$\l^3$  is a real
manifold of dimension $3n+4.$ The convergence in $\Mg_n$ means the one of {\em marked} conformal structures in ${\cal T}_n$ and Weierstrass data. The ordered sequence of points in the 
mark and the logarithmic growth at infinity provide global analytical coordinates with the same underlying topology. This space is a
$(n+1)!$-sheeted covering of the space $\Gg_n$ of (non marked) entire maximal 
graphs over $\{x_3=0\}$ with vertical limit normal vector at the end and $n+1$ conelike singularities. The underlying topology of $\Gg_n$ is equivalent to the uniform convergence of graphs on compact subsets of $\{x_3=0\}.$
Moreover, identifying marked entire maximal 
graphs differing by  translations,  rotations
about  a vertical axis, homotheties or symmetries about a horizontal plane preserving the marks, the
quotient space $\hat{\Mg}_n$ is an analytic manifold of dimension $3n-1$
diffeomorphic to the spinorial bundle ${\cal S}_n.$}
\end{quote}
{\bf Acknowledgments:} {\small We would like to thanks to Antonio Ros for helpful conversations during the preparation of this work. We are also indebted to Michael Wolf for some useful comments. 
This paper was carried out during the third author's visit at the Departamento de Geometría y Topología de la Universidad de Granada (Spain), in February-April, 2003. The third author is grateful to the people at the department for their hospitality.}

\section{Notations and Preliminary results} \label{sec:prelim}
Throughout this paper, $\l^3$ will denote the three dimensional
Lorentz-Minkowski space $(\r^3,\langle , \rangle),$ where $\langle ,
\rangle=dx_1^2+dx_2^2-dx_3^2.$
We say that a vector ${\bf v} \in \r^3- \{ {\bf 0} \}$ is spacelike,
timelike or lightlike if
$\langle {\bf v}, {\bf v} \rangle$ is positive, negative or zero,
respectively. The vector
${\bf 0}$ is spacelike by definition.  A plane in $\l^3$ is spacelike,
timelike or lightlike if the induced metric is Riemannian, non degenerate
and indefinite or degenerate, respectively.
Throughout this paper, $\overline{\c}$ denotes the extended complex plane
$\c \cup \{\infty\}.$
We denote by $\h^2 = \{ (x_1,x_2,x_3) \in \r^3 \;:\;
x_1^2+x_2^2-x_3^2=-1\}$ the hyperbolic sphere in $\l^3$ of constant
intrinsic curvature $-1.$ Note that $\h^2$ has two connected  components,
one on which $x_3 \geq 1$ and one on which $x_3 \leq -1.$ The stereographic
projection $\sigma$ for $\h^2$ is defined as follows:
$$\sigma:\overline{\c} - \{|z|=1\} \longrightarrow \h^2 \,; \; z \rightarrow
\left(\frac{2 \mbox{Im} (z)}{|z|^2-1}, \frac{2 \mbox{Re} (z)}{|z|^2-1},
\frac{|z|^2+1}{|z|^2-1} \right),$$ where $\sigma(\infty)=(0,0,1).$

An immersion $X:M \longrightarrow \l^3$ is spacelike if
the tangent plane at any point is spacelike.  The Gauss map $N$ of $X$
(locally well defined) assigns to each point of $M$ a point
of $\h^2 .$ If $X$ is spacelike, $N$ is globally well defined (that is to
say, $M$ is orientable) and $N(M)$ lies in one of the components of $\h^2.$
A maximal immersion $X:M \longrightarrow \l^3$ is a spacelike immersion
such that its mean curvature vanishes. Using isothermal parameters
compatible with a fixed orientation $N:M \to \h^2,$ $M$ has in a natural
way  a conformal structure, and
the map $g \df \sigma^{-1} \circ N$ is meromorphic.
Moreover, there exists a holomorphic $1$-form $\phi_3$ on $M$ such that the
$1$-forms
$\phi_1=\frac{i}{2} \phi_3(\frac{1}{g}-g),$ $\phi_2=-\frac{1}{2} \phi_3 (\frac{1}{g}+g)$
are holomorphic, and together with $\phi_3,$ have
no real periods on $M$ and no common zeroes. Up to a translation, the immersion is given by 
$$X= \mbox{Re}  \int ( \phi_1,\phi_2,\phi_3).$$
The induced Riemannian metric $ds^2$ on $M$ is given by
$ds^2=  |\phi_{1}|^2 +|\phi_{2}|^2- |\phi_{3}|^2 =\left( \frac{|\phi_3|}{2}
(\frac{1}{|g|}-|g|) \right)^2.$
Since $M$ is spacelike, then $|g| \neq 1$ on $M.$ 

\begin{remark}
For convenience, we also deal with surfaces $M$ having $\partial (M) \neq \emptyset,$ and in this case, we always suppose that $\phi_3$ and $g$ extend analitically beyond $\partial M.$
\end{remark}
Conversely, let $M,$ $g$ and $\phi_3$ be a Riemann surface with possibly non
empty boundary, a meromorphic map on $M$ and an holomorphic $1$-form $\phi_3$
on $M,$  such  that $|g(P)| \neq 1,$ $\forall P \in M,$ and the $1$-forms
$\phi_j,$ $j=1,2,3$ defined as above are holomorphic, have no real periods
and have no common zeroes. Then the conformal immersion $X= \mbox{Re}  \int
( \phi_1,\phi_2,\phi_3)$  is maximal, and its Gauss map is
$\sigma \circ g.$ We call $(M, \phi_1,\phi_2,\phi_3)$ (or simply
$(M,g,\phi_3)$)  the Weierstrass representation of $X$.
For more details see, for instance, \cite{kobayashi}.
If we allow $A_X:=\{p \in M \;:\; |g(P)|=1\} \cup \{P \in M
\;:\;\sum_{j=1}^3 |\phi_j|^2(P) = 0 \} \neq \emptyset,$ we say that $X:M
\rightarrow \l^3$  has singularities,  and that  $A_X$ ( resp., $X(A_X)$) is the
set of singularities of $X$ (resp., $X(M)$). Observe that  $ds^2|{A_X}=0,$ and in
particular, the map $X$ is not an immersion at the singular set $A_X.$
A maximal  immersion $X:M \to \l^3$ with a compact set of
singularities $A_X \subset M$  is said to be complete if any divergent path
in $M$ has infinite length.

\begin{remark} \label{re:cambio}
The transformation $(M, \Phi_1,\Phi_2,\Phi_3) \rightarrow (M, i \Phi_1,i
\Phi_2,\Phi_3)$ converts Weierstrass
data of maximal  surfaces in $\l^3$ into Weierstrass data of
minimal surfaces in $\r^3,$ and vice versa, provided that the period
problem is solved in each case. For more details about minimal surfaces,
see \cite{osserman}.
\end{remark}

\subsection{Behaviour of maximal  surfaces around isolated
singularities.}

In this subsection we study some basic properties of maximal surfaces around isolated singularities.
Our analysis includes the general immersed case, emphasizing the topological behaviour of singular points.

Let $X: {\cal D}\to \l^3$ be a continuous map  defined on an open  disk
${\cal D},$ $q$ a point in ${\cal D},$ and suppose $X$ is a maximal
immersion on ${\cal D}-\{q\}.$ Let $z$ be a conformal parameter
on ${\cal D}-\{q\}$ associated to the metric $ds^2$ induced by $X,$ and
write $ds^2= h(z) |dz|^2,$ where $h(z)>0$ for any  $z \in z({\cal
D}-\{q\}).$  By definition, $q$ is an {\em (isolated) singularity} of $X$
if for any  sequence $\{q_n\} \subset {\cal D} -\{q\}$ tending to $q,$  the
limit $\lim_{n \to \infty} h(z(q_n))$ vanishes. In this case, we say that
$X({\cal D})$ is a maximal  surface with a singularity at $X(q).$

There are two kinds of isolated singularities: {\em branch points} and {\em
special singularities}.
In case ${\cal D}-\{q\}$ endowed with the induced complex structure is
conformally a once punctured
disc, then $q$ (respectively, $X(q)$) is said to be a {\em branch point} of
$X$ (respectively, of $X({\cal D})$). This  means that, although the metric
degenerates at $q,$ the Weierstrass data $(g, \phi_3)$ of $X,$ and so the
Gauss map, extend meromorphically to $q,$  see for instance \cite{Es}.
In this case, $\phi_1,$ $\phi_2$ and $\phi_3$ vanish at $q$ (that is to
say, $A_X=\{q\}$) and $X$ is not an embedding around the singularity (see
Remark
\ref{re:branch}). The local behavior at the singularity is similar to
the case of minimal surfaces in $\r^3$ (see \cite{osserman}) ).

To avoid false branch points or trivial coverings without geometrical
significance, we always assume that the set of self intersections of $X$
is either empty or a one dimensional analytic variety. In other words, we
suppose that the immersion can not be factorized by a nontrivial covering.

Suppose now that ${\cal D}-\{q\}$ is  conformal to an annulus $\{z \in \c
\;:\; 0<r < |z| < 1\}.$ In this case $X$ extends continously to $C_0=\{z
\in \c \;:\; 0<r < |z| \leq 1\},$ with $X(\{|z|=1\}) = X(q).$ If
$J(z)=1/\overline{z}$
denotes the inversion about $\{|z|=1\},$ then by Schwarz reflection $X$
extends analytically to
$C=\{z \in \c \;:\; r < |z| < 1/r \}$ and satisfies $X\circ J = -X +
2P_0,$ where $P_0=X(q).$ Labelling $(g,\phi_3)$ as the Weierstrass
representation
of the extended immersion, it follows straightforwardly that
$J^*(\phi_k)=-\overline{\phi_k}, \,\,k=1,2,3,$ and therefore that
$g \circ J=1/\overline{g}$ on $C.$ In particular $|g|=1$ on $\{|z|=1\},$
that is to say, $A_X=\{|z|=1\}.$
Moreover, observe that $dg \neq 0 \;\;\mbox{on} \;\;\{|z|=1\}.$ Indeed, the
critical points of $g$ on the set $|g|^{-1}(1)$ correspond to the cross
points of the
nodal set of the harmonic function $\log|g|,$ but there are no nodal curves
of this function in $C_0-\{|z|=1\}$ since the surface is spacelike on this
domain.
In this case, the Gauss map of $X$ has no well defined limit at $q,$ and
$q$ (respectively, $P_0$) is said to be an {\em special singularity} of
the immersion $X$ (respectively, of $X({\cal D})).$ As we will see during the
proof of the following lemma, the local behavior at an special
singularity is
determined by two integers, namely the degree $m$ of the map $g:
\{|z|=1\}\to \{|z|=1\},$ which we will call the {\em degree} of $g$ at the
singularity, and
the number $n$ of zeros of
$\phi_3$ on $\{|z|=1\},$ to which we will refer as the {\em vanishing
order} at the singularity. 
It is interesting to notice that the vanishing order $n$ is always even. Indeed, 
the zeros of
$\phi_3$ on $\{|z|=1\}$ are the critical points of the harmonic function
$x_3$ on this set
and therefore correspond to the cross points of the nodal set of this function
on $\{|z|=1\}$ (up to translation, we are assuming that $P_0$ coincides
with the origin). Furthermore, recall that the multiplicity of a point
$z_0$ as zero of $\phi_3$ is equal to the number of nodal curves meeting at
$z_0$ minus one. By the maximum
principle there are no domains in $C_0$ bounded by nodal curves, and $x_3$
changes
sign when crossing a nodal curve. These facts show that there is an even number
of nodal
curves crossing the circle  $\{|z|=1\},$ i.e., $n = 2k.$
Examples of maximal  discs with an special singularity are given
by the following Weierstrass data on $C_0:$ $$g=z^m,\;\; \phi_3=i
\frac{(z^{2
k}-1)}{z^{k+1}}dz,$$ $m \geq 1,$ $k\geq 1.$  Call $X:C_0 \to \l^3$ the
associated conformal immersion and $D_{m,k}=X(C_0).$
When $k=0,$ we take $g=z^m,$ $\phi_3=\frac{dz}{z}$ and define  $D_m \equiv D_{m,k}$
as above.

For instance, in the later case, $m=1$ corresponds to the Lorentzian catenoid. 
The asymptotic behavior of maximal surfaces around special singular points with vanishing number $n=0$ has been extensively studied in \cite{kly-mik}. In the following lemma we fix attention on some topological properties of general isolated singularities, and express them in terms of the Weierstrass representation. 

In the case of embedded singularities, we emphasize that the degree of $g$ and the vanishing order at the singularity satisfy: $m=1,$ $k=0.$ We also include, just for completeness, a brief reference to the asymptotic behavior of embedded singularities (conelike singularities), although a more extensive study of these analytical properties can be found in \cite{kobayashi}, \cite{ecker} and \cite{kly-mik}.

\begin{lemma}[local structure of special singularities] \label{lem:sing}
Let  $X: {\cal D}\to \l^3$ be a maximal  immersion defined on an
open  disk ${\cal D},$
and suppose that $X$ has an special singularity at $q \in {\cal D}.$ Let
$X:\{0<r<|z|\leq 1\} \to \l^3$ be a conformal reparameterization of $X$ with
$P_0:=X(q)=X(\{|z|=1\})$ the singular point. As above denote by $m$ and $2k$ the
degree of $g$ and the vanishing order at the singularity, respectively.
Let $\Pi$ be a spacelike plane containing $P_0$ and label $\pi: \l^3 \to
\Pi$ as
the Lorentzian orthogonal projection. 
Then, there exists a small closed
disc $U$ in $\Pi$  centered at $P_0$  such that $(\pi \circ X):V-\{|z|=1\}
\to
U-\{P_0\}$ is a covering of $m+k$ sheets, where $V$ is the connected
component of $(\pi \circ X)^{-1}(U)$ containing $\{|z|=1\}.$
In particular, there exists $r_0 \in [r,1[$ such that  $X:\{0<r_0<|z|<1\}
\to \l^3$ is an embedding if and only if $g$ is injective on $\{|z|=1\}$
and $\phi_3(z)\neq 0,$ $|z|=1,$ (that is to say, $m=1,$ $k=0$). In this
case,  $X(\{0<r_0<|z| \leq 1\})$ is a graph over $\Pi$ contained in one of the halfspaces determined by this plane and  asymptotic to the light cone at the singularity.
The point $P_0$ is then said to be a conelike singularity.
\end{lemma}

\begin{figure}[htbp] 
\begin{center}
\includegraphics[width=10cm]{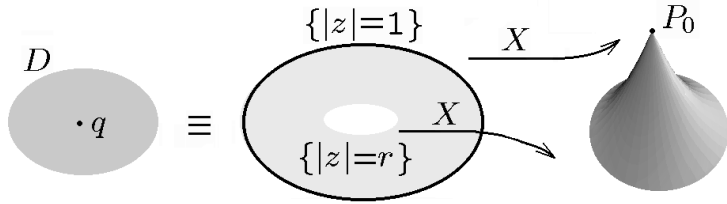}
\caption{A conelike singularity. \label{fig:conelike}}
\end{center}
\end{figure}

\begin{proof}
Up to an ambient isometry we can assume that $\Pi = \{x_3 = 0\}$ and $P_0$
coincides with the origin. In the sequel, we put $D=X(\{0<r<|z| \leq 1\}).$
Write $g(z)=w(z)^m$ in a small neighborhood of $\{|z|=1\},$ where $w(z)$ is holomorphic. Taking $w$ as a
new parameter, $g(w)=w^m$ and $\phi_3(w)$ are well defined in a open
annulus $C_1$ in the $w$-plane containing $\{|w|=1\}.$ Fix $C_2=\{0<r_1<|w|<1/r_1\}
\subset C_1,$ $r_1$ small enough and put $C_3=\{0<r_1<|w|<1\}.$ If $k \neq
0,$ label $a_1,\ldots, a_{2 k}$  as the zeros of
$\phi_3$ on $\{|w|=1\},$ where each zero appears as many times as its
multiplicity. Then, we can write
$$\phi_3(w)= i \frac{\prod_{j=1}^{2 k} (w-a_j)}{w^{k+1}} f(w) dw,$$ for a
suitable non-vanishing holomorphic function $f(w).$ For $k=0$ we simply
write: $\phi_3(w)=  \frac{f(w)}{w} dw.$
\vspace{0.2cm}

\begin{quote}{ {\bf Claim:} \em There exists a small closed disc $U$ in
$\{x_3=0\}$ centered at the origin such that the connected component $V$ of
$(\pi \circ X)^{-1}(U)   $ containing $\{|w|=1\}$ satisfies $X(\partial
V -\{|w|=1\})\subset \pi^{-1}(\partial U).$} \end{quote}

Indeed, assume by contradiction that there exists a positive sequence
$\{s_h\} \to 0$  such that, for each $h,$ $X(\partial V_h)$ is not
contained in
$\pi^{-1}(\partial U_h),$ where $U_h$ is the disc of radius $s_h$
centered at the origin in $\{x_3=0\}$  and  $V_h$ is the connected
component of
$X^{-1}(D \cap \pi^{-1}(U_{h}))$ containing $\{|w|=1\}.$

Then, $V_h$ is a connected domain containing $\{|w|=1\}$ and a piece of
$\{|w|=r_1\}$ in its boundary. Then, the set 
$X^{-1}(D \cap
\{(0,0,x) \;:\; x \in \r\})$  contains  a curve joining the two
components of $\partial C_3.$ This contradicts that $X$ is
spacelike in
$C_3-\{|w|=1\}$ and proves the claim.

\vspace{0.2cm}

Let us now show that $\pi \circ X:V-(\pi \circ X)^{-1}(0) \to U-\{0\}$ is a
finite covering. As $X$ is spacelike on $C_3-\{|w|=1\},$ then $\pi \circ X$
is a local diffeomorphism on this set. Moreover, since $(\pi \circ
X)(\partial V)$ is contained in $\partial U,$ it follows that
$\pi \circ X:V-(\pi\circ X)^{-1}(0) \to U-\{0\}$ is a local diffeomorphism.
As $(\pi\circ X)^{-1}(0)\cap V$ is compact, it is clear that $\pi \circ
X:V-(\pi \circ X)^{-1}(0) \to U-\{0\}$ is proper and so
we infer that $\pi \circ X:V-(\pi \circ X)^{-1}(0) \to U-\{0\}$ is a finite
covering. Since $U-\{0\}$ is a cylinder,
we deduce that $V-(\pi \circ X)^{-1}(0)$ is a cylinder too and $(\pi \circ X)^{-1}(0)\cap V=\{|w|=1\}.$

Now observe that, by elementary topology, the number of sheets of the
covering $\pi \circ X: V-(\pi \circ X)^{-1}(0) \to U-\{0\}$ is equal to
the
winding number around the origin of the planar curve $ (\pi \circ X)
(w(\theta)),$ where
$w: \theta \in [0,1] \to \partial V$ is a parameterization of
the loop $\gamma:= \partial V - \{|w|=1|\}.$ Since $ \pi \circ X (\gamma)$
is the round circle $\partial U,$
this winding number coincides, up to dividing by $2\pi$,  with the
variation of the argument of the complex function  $\theta \to
\frac{d}{d\theta}
(\pi \circ X) (w(\theta)).$
We have:
$$\frac{\partial}{\partial \theta}(x_1 + i x_2)(X(w(\theta)) =
-i \frac{\overline{F(w)w'}}{2\overline {w}^m} \left( 1 +
\frac{F(w)w'}{\overline{F(w)w'}} |w|^{2m}\right)(\theta),$$
where $F(w)= i \frac{\prod_{j=1}^{2 k} (w-a_j)}{w^{k+1}} f(w) $ if $k \geq
1$ and $F(w)= \frac{f(w)}{w}$ if $k=0.$
Since $|w(\theta)| < 1,$ it follows that $\mbox{Real}\left( 1 +
\frac{F(w)w'}{\overline{F(w)w'}} |w|^{2m}\right)(\theta) > 0,$ and
$\mbox{Real}(\frac{w}{a_j}-1)(\theta) < 0,\, j=1,\ldots,2k.$ Moreover,
as $J^* (\phi_3) = - \overline{\phi_3},$
then $\overline{f(1/\overline{w})}= - (\Pi_{j=1}^{2k} a_j) f(w)$  for $k\ge
1$ and $\overline{f(1/\overline{w})}=  f(w)$ if $k=0$.  Therefore the
argument of
$f$ is constant on $\{|w|=1|\}.$ Since in addition $f$ is holomorphic and
non-vanishing in $C_3,$ we infer that the variation of its argument is also
zero on
$\partial V.$
It follows that the variation of the angle under consideration  on $[0,1]$  is
equal to $m+k.$ This proves the first part of the lemma.
For the second one, note that $X:V-(\pi \circ X)^{-1}(0) \to \l^3$ is
an embedding if and only if the continuous map $x_3$ separates the fibers
of the
covering $\pi \circ X:V-(\pi \circ X)^{-1}(0) \to U-\{0\}.$ This is
equivalent to saying that this covering is a homeomorphism, that is to say, $m=1$ and
$k=0.$ Therefore, ${\cal D}$ is a graph over $\{x_3=0\}$ locally around the singular point.
Since $k=0,$ we know that there are no
interior zeros of $x_3$ in $C_3$ close to $\{|w|=1|\},$ and so $D$  lies above or below $\Pi.$
Recall that
$(g,\phi_3)=(w,f(w) dw/w),$ where $f(w)$ is real and non vanishing on
$\{|w|=1\}.$ Writing the immersion in polar coordinates $w=r e^{i \theta}
,$ it is
straightforward to check that:
$$2 X(r e^{i \theta})=\mbox{Real} \int_{1}^r \frac{f(s e^{i \theta})}{s}
\left(i (\frac{e^{-i \theta}}{s}-s e^{i \theta}),-(\frac{e^{-i
\theta}}{s}+s e^{I \theta}),
2 \right).$$ An elementary computation gives:
$$\lim_{r \to 1} ||\frac{X(r e^{i \theta}) }{x_3(X(r e^{i
\theta}))}-(\sin(\theta),-\cos(\theta),1)||_1=0,$$ 
where $||\cdot||_1$ is the ${\cal C}^1$ norm in ${\cal C}^1([0,2 \pi],\r^3),$ which proves that $D$ is asymptotic to the half light cone at the singularity. This concludes the proof.
\end{proof}

\begin{remark}[local structure of branch points] \label{re:branch}
Let $X: \{ |z| \leq 1\} \to \l^3$ be a conformal maximal 
immersion branched at the origin, and call $(g,\phi_3)$ its
Weierstrass data. Up to a Lorentzian isometry we can suppose $|g|<1$ and
$\frac{\phi_3}{g}
= z^k h(z)dz,$
where $k\geq 1$ and $h(0) \neq 0.$ With the same notations as in Lemma
\ref{lem:sing},  the map $\pi \circ X$ is a covering of degree $k+1$
branched at the origin. In particular $X$ is never
an embedding. The proof of this fact
uses the same arguments as in the previous lemma.
\end{remark}

\subsection{Global behavior of complete embedded maximal  surfaces
with isolated singularities} \label{sec:global}
Let  $M$ be a differentiable surface without boundary, $X: M \to \l^3$  a
continuous map and $F=\{q_1,\ldots,q_n\} \subset M$ a finite set. We say
that $X$ is
a complete  maximal  immersion with singularities at the points
$q_1,\ldots,q_n,$ if
$X: M-F\to \l^3$ is a maximal  immersion, the points
$q_1,\ldots,q_n,$ are isolated singularities of $X,$ and   every divergent
curve in $M$ has infinite length for the induced (singular) metric.

\begin{proposition}  \label{pro:graph}
Let $X: M \to \l^3$ be a maximal  immersion with a finite set of
singularities $F,$ and let $\Pi$ be any spacelike plane in $\l^3.$
Then the following statements are equivalent:
\begin{enumerate}
\item $X$ is an embedding and complete,
\item $X$ is complete and all its singularities are of conelike type,
\item $X(M)$ is an entire graph over $\Pi.$
\end{enumerate}
In this case, if $C$ denotes a closed disc in $M$ containing $F,$ then
$M-\stackrel{\circ}{C}$ is conformally equivalent to a compact Riemann
surface (with non empty boundary) minus an interior point which we call the
end of $M.$ Moreover, the Weierstrass data $(g,\phi_3)$ of $X$ extend
meromorphically to the end of $M,$ and so, the Gauss map and the tangent
plane are well defined at the end.
\end{proposition}

\begin{proof}
Label $\pi:\l^3 \to \Pi$ as the Lorentzian orthogonal projection over $\Pi.$
Up to an ambient isometry, we can assume that $\Pi=\{x_3=0\}.$
To check that $(1)$ implies $(2)$ just see Lemma \ref{lem:sing} and Remark
\ref{re:branch}.
Suppose $X$ is complete and its singularities are of conelike type.
Since $X$ is spacelike out of the singular points,
then $\pi \circ X:M-F \to \{x_3=0\}$ is a local diffeomorphism.
Call $ds^2$  the metric induced on $M$
by $X$ and $ds_0^2$ the Euclidean metric on $\{x_3=0\}.$ Then $ds^2 \leq
(\pi\circ X)^*(ds_0^2)$ outside the singularities. By Lemma \ref{lem:sing} and the
completeness of  $X,$  the map $\pi \circ X$ has the
path-lifting property. Hence $\pi \circ X$ is an unbranched covering of the plane,  and
therefore, a homeomorphism. This proves that  $(2)$ implies $(3).$

Assume now that $X(M)$ is a graph over $\Pi,$ then it is clear that $X$ is
an embedding. By Lemma \ref{lem:sing} and Remark \ref{re:branch} it is
straightforward that the singularities of $X$ are of conelike type.
Without loss of generality we assume that $|g|<1$ on $M-F,$ because
$X|_{M-F}$ is spacelike.
On the other hand, observe that $(\pi\circ X)^*(ds_0^2)\leq ds_1^2 \leq
|\frac{\phi_3}{g}|^2,$ where
$ds_1^2 = |\phi_1|^2 + |\phi_2|^2 + |\phi_3|^2 = \frac{1}{4} |\phi_3|^2
(|g| + 1/|g|)^2.$
Thus, the flat metric $|\frac{\phi_3}{g}|^2 $ is complete, and so it
follows from  classical
results of Huber \cite{hub} and  Osserman \cite{osserman} that
$M-\stackrel{\circ}{C}$  is conformally a once punctured compact Riemann surface with compact boundary, $\phi_3/g$ has  poles at the puncture
and $g$ extends holomorphically to the puncture.
Consequently, $|g|<1-\epsilon$ on $M-\stackrel{\circ}{C}$  for some
$\epsilon >0.$ Therefore, on $M-\stackrel{\circ}{C},$
the induced metric $ds^2 = \frac{1}{4} |\phi_3|^2 (|g| - 1/|g|)^2 \geq
\frac{\epsilon^2}{4} |\frac{\phi_3}{g}|^2$
and so it is complete.
\end{proof}

Following  Osserman \cite{osserman}, the previous proposition shows that
any complete and embedded maximal  surface with a finite set of
singularities has finite total curvature outside any neighborhood of the
singularities.

In the sequel and for the sake of brevity we write \textbf{CMF} to refer
to a complete maximal 
immersion (surface, graph,$\ldots$) with a finite set of special
singularities.
\begin{remark}
We emphasize here that there exist entire spacelike graphs over
$\{x_3=0\}$ in $\l^3$ which are not complete for the induced metric.
The previous proposition asserts that, in the maximal case, an entire graph
is always complete.
\end{remark}
The next lemma describes the asymptotic behaviour of  a CMF graph. The proof follows from
some well known classic results by Osserman \cite{osserman}, Jorge and Meeks \cite{jorgemeeks} and R. Schoen \cite {schoen} for minimal surfaces. A different approach in
the Lorentzian setting can be found in \cite{kly}. We omit the proof.
\begin{lemma} \label{lem:graph}
Let $X:M \rightarrow \l^3$ be a maximal embedding 
with a finite set $F$ of singularities and vertical limit normal vector at the end,
and label $(g,\phi_3)$ its Weierstrass representation.

Then, $\mbox{Max}\{O_\infty (\phi_k), \, k=1,2,3\}=2,$ where $O_\infty (\phi_k)$ is the pole
order of $\phi_k$ at the end, $k=1,2,3.$
 
Moreover, if we write $G:=X(M)=\{u(x_1,x_2) \;:\; (x_1,x_2) \in \r^2\}$ as an entire graph over $\{x_3=0\},$ then:
$$u(x_1,x_2)= c \log{|(x_1,x_2)|}+b +\frac{a_1 x_1+a_2 x_2}{|(x_1,x_2)|^2}+
O(|(x_1,x_2)|^{-2}),$$ for suitable constants $c,$ $b$, $a_1$ and $a_2.$
\end{lemma}
Take  an arbitrary complete embedded maximal  surface with a
finite set of singularities, and up to a Lorentzian isometry, suppose that
the limit
normal vector at its unique end is vertical. From Lemma \ref{lem:graph}, the surface is an entire graph $u(x_1,x_2)$ over
$\{x_3=0\}$ which is asymptotic either to a vertical half catenoid ({\em
catenoidal end}) or to a horizontal plane ({\em planar end}). The
asymptotic behavior of the surface is controlled by the real constant $c$
appearing in Lemma \ref{lem:graph}, called the {\em logarithmic
growth}
of the end (or of the surface).  Note that $c=0$
(respectively, $c \neq 0$) if and only if the end is of planar type
(respectively, of catenoidal type).
\vspace{0.2cm}

Let $X:M \to \l^3$ be a maximal  surface with a finite set $F$ of
singularities, and fix an orientation on $M.$
\begin{definition}[Flux and torque of closed curves]
Let $\gamma(s):[0,L] \to M-F$ be an {\em oriented} closed curve
parameterized by arclength. Denote by $\nu$ the unit conormal to the curve
such that $(\nu, \gamma')$ is positively oriented with respect to the
orientation of $M.$
The flux of $X$ along $\gamma$ is defined by:
$$F(\gamma)= \int_{\gamma} \nu(s) ds,$$
and the torque of $X$ along $\gamma$ is defined by:
$$T(\gamma)= \int_{\gamma} X(\gamma (s))\wedge \nu (s) ds,$$
where $\wedge$ refers to the Lorentzian exterior product.
\end{definition}
Since $X$ is harmonic, it follows from Stokes theorem that $F(\gamma)$
depends only on the homology class of $\gamma$ on $M-F.$ Also, observe that
since the immersion is maximal then the
divergence of the vector valued 1-form defined on $M$ by: $\alpha_p(v) =
X(p)\wedge ^\ast dX_p (v)$, where $p\in M$ and $v \in T_p M,$ is zero. Therefore
the torque depends only on the homology class in $M-F$ of the curve.
However, note that the torque depends on the choice of the origin in $\l^3.$

Suppose now that $X: M \to \l^3$ defines a CMF graph with vertical limit normal vector at the end, and call
$\{q_1,\ldots,q_n\}\subset M$ its
singularities. Up to
an ambient isometry we can assume that the end is asymptotic either to a
vertical upward half-catenoid or to a horizontal plane.
Henceforth, we will orient $X$ so that the  normal at the end point
downwards, that is to say the Gauss
mapping satisfies $|g| < 1$.
Consider a compact domain $\Omega \subset M-F$ such that
$\partial \Omega = \gamma_{\infty}\cup\gamma_1\cup\ldots \cup\gamma_n$ with
$\gamma_i$
being the boundary of a closed disc in $M$ containing $q_i$ in its interior
and no other singularity and $\gamma_{\infty}$ the boundary of a closed
disc containing all the singularities in its interior.
The orientation on $\Omega$ induces an orientation on these curves in the
usual way, that is to say,
$(\nu_i,\gamma_i')$ is positively oriented where $\nu_i$ is the exterior
unit conormal to $\gamma_i.$
Call $F_{\infty}$  and $F_i$ (resp., $T_{\infty}$ and $T_i$) the flux  (resp., torque) of
$X$ along these oriented curves $\gamma_{\infty}$ and
$\gamma_i,$ respectively, $i=1,\ldots,n.$
Note that 
$F_\infty=0$ for planar ends, and $F_\infty=2 \pi (0,0,c)$ for a
catenoidal end with logarithmic growth $c.$ Moreover, if we take a sequence of closed curves $\{\gamma_i^k\}_{k \in \n}$ in the same homology class of $\gamma_i$ tending, as $k \to \infty,$  to the singularity $q_i,$ it is not hard to check that the conormal directions  along this curves go to the same component of the light cone. Therefore, it is easy to check that the flux $F_i$ at $q_i$ is a timelike vector, $i=1,\ldots,n.$

The geometrical interpretation of $T_i$ and $T_\infty$ is similar to the one in the case of minimal surfaces. For instance, it is not difficult to prove that $T_i= X(q_i) \wedge F_i,$ for any $i.$ Moreover,  if the end is of catenoidal type, then $T_\infty=P_\infty \wedge F_\infty,$ where $P_\infty$ is any point on the axis of the catenoidal end. If the end is planar and not of Riemann type, $T_\infty=0,$ and in the case of planar  ends of Riemann type,  $T_\infty$ is a horizontal vector contained in the asymptotic line of the end. We omit the proof of all these details and refer to  \cite{Perez} and \cite{Kusner} for the analogous results in the theory of minimal surfaces. These references are also of interest when we deal with the following equilibrium relations:

\begin{proposition} \label{pro:equi}
With the above notations:
\begin{enumerate}[(1)]
\item $F_{\infty} + \sum_{i=1}^{n} F_i = 0,$
\item $T_{\infty} + \sum_{i=1}^{n} T_i = 0.$
\end{enumerate}
\end{proposition}
\begin{proof}
(1) follows from harmonicity of $X$ and Stokes theorem. To prove (2), take into account the maximality and apply
the divergence theorem to the vector valued 1-form defined on $M$:
$\alpha_p(v) = X(p)\wedge ^\ast dX_p(v)$, $p\in M,$ $v \in T_p M.$
\end{proof}
Formula (1) gives the following:
\begin{corollary} \label{co:halfspace}
There does not exist any CMF  embedded surface having a planar end and all
the singularities
pointing upwards (reps. downwards). Likewise  there is no CMF embedded
surface
having an upward (downward) vertical catenoidal end and all the
singularities pointing upward (resp. downward).
\end{corollary}
We can now give the following characterization of the Lorentzian catenoid.

\begin{theorem}
The Lorentzian half catenoid is the only CMF embedded surface having
downward pointing singularities all with vertical flux.
\end{theorem}

\begin{proof}
Denote by $X: M \to \l^3$ such an immersion, $q_1,\ldots,q_n \in M$ its
singularities.
From Proposition \ref{pro:equi} and Corollary \ref{co:halfspace},  and up
to a ambient isometry, the end of $X$ must be a vertical upward half
catenoid. Conformally we
can put $M-\{q_1,\ldots,q_n \}= \Sigma -(\partial \Sigma \cup
\{q_{\infty}\}),$ where
$\Sigma$ is a compact planar domain with analytic boundary, and
$q_{\infty}$ is an interior
point of $\Sigma$ (which corresponds to the end). It is clear that
$\partial\Sigma$ consists
of $n$ analytic circles $a_1,\ldots,a_n$ which correspond to
the singularities. Label $(g,\phi_3)$ as the Weierstrass data of $X: \Sigma
-\{q_{\infty}\}\to \l^3$ and observe that the hypothesis
imply that $\phi_1$ and $\phi_2$ are exact. It follows that the conformal
minimal immersion $Y: \Sigma -\{q_{\infty}\}\to \r^3$ given by
$Y=\mbox{Real}\int (i\phi_1,i\phi_2,\phi_3)$
is well defined. Since along $a_i$ the Gauss map is horizontal and
injective and $x_3\circ Y$ is constant, then  the curve $Y(a_i)$ is a
convex planar geodesic,  for any $i$.
Furthermore the end of the immersion $Y$ is an upward vertical catenoidal
end. Call $D_i$
the planar closed disc bounded by $Y(a_i)$ and consider the
topological surface $S$
obtained by attaching  $D_i$ to $\Sigma$ (we  identify any point in $a_i$ and
its image under $Y$ in $Y(a_i)\subset
D_i$). Let us consider the continuous map $Z:S\to \r^3$ given by
$Z|_{\Sigma}=Y$ and
$Z|_{D_i}= \mbox{id},$ for any $i.$ If $\pi: \r^3\to \{x_3 =0\}$ denotes
the orthogonal
projection, it is not hard to see that $\pi \circ Z: S\to \{x_3 =0\}$ is a
local homeomorphism (recall that, without loss of generality, we can assume
that $|g|<1$ on $\Sigma-\partial \Sigma$).
Moreover since the end of $Y$ is a vertical catenoid, $\pi\circ Z$ is
proper. Thus $\pi \circ Z$
is a homeomorphism. In particular $Y(\Sigma)$ is a graph over $\{x_3 =0\}.$

Then, $Y:\Sigma \to \r^3$ is a minimal embedding with a vertical upward
catenoidal end, $Y(\partial \Sigma)$ consisting of a finite set of
horizontal planar convex geodesics with vertical  downward fluxes. Under
these assumptions, Ros \cite{ros} has proved that $Y(\Sigma)$ is
necessarily a half catenoid.
\end{proof}

\subsection{ Existence and Uniqueness of CMF graphs with any number of
singularities}\label{sec:exist}
Although we have stated it only in the three dimensional case, the following existence and uniqueness result  is also valid for maximal graphs in $\l^{n+1}.$
\begin{theorem}[Klyachin \cite{kly}] \label{th:uniqueness}
Let $h(x)$ be a function defined on a compact set $K \subset \{x_{3}=0\},$ $\mbox{Int}(K)=\emptyset,$  and satisfying the inequality $|h(x)-h(y)| <
|x-y|$ for all $x,$ $y \in K,$ $x \neq  y.$ Then, for every timelike vector $v$  there exists a unique
solution $u \in C^2(\r^2 -K) \cap C(\r^2)$ to the  maximal graph equation such that $u|_K =h$ and $F_\infty=v,$ where as above $F_\infty$ is the flux at
the end.
\end{theorem}

If $G=\{(x,u(x)) \, x \in \r^2\}$ is an entire maximal graph  with a finite number of singularities, it is not hard to check that  $$|u(x)-u(y)| <
|x-y|, \quad \mbox{for all}\; x,\,y \in \r^2, \, x \neq  y.$$  
Therefore, the previous theorem implies that any  CMF graph $G$ with vertical limit normal vector at infinity is uniquely determined by the position of its singular points  and  its logarithmic growth at infinity.
\begin{remark} \label{re:uniqueness}
As a consequence of the previous theorem, the group of  ambient isometries preserving an embedded CMF
surface coincides with:
\begin{itemize}
\item
the group of ambient isometries leaving the set of its singularities
invariant and preserving a halfspace containing the surface in case the
surface has a catenoidal end,
\item
the group of ambient isometries leaving the set of its singularities
invariant  in case the surface has a planar end.
\end{itemize}
\end{remark}
Notice that the existente result in \cite{kly} is very general but implicit. In particular, one can not control whether the  isolated singular points are upward pointing, downward pointing or regular. Theorem \ref{th:exist} below gives an explicit method for constructing CMF graphs which allows us to control the geometry of the arising examples.  Moreover, it shows that  CMF graphs and meromorphic data on compact Riemann surfaces admitting a mirror involution are closely related. This fact will be crucial for understanding the moduli space of this kind of surfaces. We start introducing the following notation. 

By definition, an open domain $\Omega \subset \overline{\c}$ is said to be
a circular domain if its boundary consists of a finite number of circles.

\begin{theorem} \label{th:exist}
Let $N$ be a compact genus $n$ Riemann surface, and let $J:N \rightarrow N$
be an antiholomorphic involution in $N.$ Assume that the fixed point set of $J$ consists of $n+1$ pairwise disjoint analytic Jordan curves $a_j,$ $j=0,1,\ldots,n,$ and that $N-\bigcup_{j=0}^n a_j=\Omega \cup J(\Omega),$ where $\Omega$ is  topologically equivalent (and so conformally) to an open planar circular domain.

Let $(g,\phi_3)$ be Weierstrass data  on $N$ such that:
\begin{enumerate}[(1)]
\item $g$ is a meromorphic function on $N$ of degree $n+1,$
$|g|<1$ on $\Omega$ and $g \circ J=\frac{1}{\overline{g}}$ on $N,$
\item $\phi_3$ is a holomorphic 1-form on $N-\{\infty,J(\infty)\},$ $\infty \in \Omega,$ with
poles of order at most two at $\infty$ and $J(\infty),$ and satisfying
$J^*(\phi_3)=-\overline{\phi_3},$
\item the zeros of $\phi_3$ in $N-\{\infty,J(\infty)\}$ coincide (with the
same multiplicity) with the zeros and poles of $g.$
\end{enumerate}
Then, the maximal  immersion $X:\overline{\Omega}-\{\infty\}
\rightarrow \l^3,$ $X(z):=\mbox{Real} \int^z \big(\phi_1,\phi_2,\phi_3
\big),$ where
$\phi_1=\frac{i}{2}(\frac{1}{g}-g)\phi_3$ and
$\phi_2=\frac{-1}{2}(\frac{1}{g}+g)\phi_3$
is well defined and $G=X(\overline{\Omega}-\{\infty\})$ is an entire maximal
graph with conelike singularities corresponding to the points
$q_j:=X(a_j),$ $j=0,$ $1,\ldots,n.$
\end{theorem}

\begin{proof}
Let us see that the map $X$ is well defined. First, notice that the curves $a_j,$ $j=0,1,\ldots,n$ generate the first
homology group on  $\overline{\Omega}-\{\infty\}.$ Moreover,
$J^*(\phi_j)=-\overline{\phi_j},$ $j=1,2,3,$ and $J$ fixes pointwise the
curves $a_j$ for any $j.$
Hence we deduce that
$$\int_{a_i} \phi_j=\int_{J(a_i)} \phi_j=
\int_{a_i} J^*(\phi_j)=-\int_{a_i} \overline{\phi_j},$$ which
means that $\phi_j$ has imaginary periods on
$\overline{\Omega}-\{\infty\},$ and so $X$ is
well defined.

On the other hand,  $|g|=1$ on $\partial \Omega$ and $(3)$ imply that  $\phi_3$ does not vanish on $a_j,$ $j=0,1,\ldots,n.$ Moreover, $\mbox{deg}(g)=n+1$ gives in addition that
$g|_{a_j}$ is injective, $j=0,1,\ldots,n,$ and so, 
by Lemma \ref{lem:sing}, all the singularities are of conelike type.
Let us prove the completeness of the metric $ds^2$ induced by $X.$ Suppose that the vanishing order of $\phi_3$ at $\infty$ and $J(\infty)$ is
$k \geq -2$ (when $k<0,$  $k=0$ or  $k>0$  this simply means that $\phi_3$ has a pole of order $-k$, is regular or has a zero of order $k,$ respectively). The classical theory of compact Riemann surfaces implies that the number
of zeros minus the number of poles of $\phi_3$ in $N$ is $2n-2$ (counting multiplicities).
Then the number of zeros of $\phi_3$ in $N-\{\infty,J(\infty)\}$ is
$2n-2-2k$ and $(3)$ implies that this is the number of poles and zeros of
$g$ in $N-\{\infty,J(\infty)\}.$ Since $\deg(g)=n+1,$ $g$ has $2(n+1)$ zeros
and poles in $N,$ and so we infer that $g$ has a zero of order $k+2$ at $\infty$
and a pole of order $k+2$ at $J(\infty)$ (take into account that  $|g|<1$ in $\Omega,$ $|g|>1$ in $J(\Omega)$ and $|g|=1$ in $\partial \Omega$). Therefore the metric $ds^2=\left( \frac{|\phi_3|}{2}
(\frac{1}{|g|}-|g|) \right)^2$ is complete.
\end{proof}
\begin{corollary}\label{co:exist}
Let $N$ be the compact genus $n$ Riemann surface:
$$N=\{(z,w) \in \overline{\c}^2 \;:\; w^2=\frac{(z-1)
\prod_{j=1}^{n} (z-c_j)}{(z+1) \prod_{j=1}^{n} (z-b_j)}\},$$ where $c_j,$
$b_j$ are pairwise distinct real numbers in $\r-\{-1,1\},$  and define on
$N$ the following meromorphic data
$$g=\frac{w-1}{w+1}, \quad \phi_3= (\frac{1}{w}-w  )dz.$$
Let $J:N \rightarrow N$ be the antiholomorphic
involution given by $J(z,w)=(\overline{z},-\overline{w}).$
Call $\{\infty_1,\infty_2\}$ the two points in $z^{-1}(\infty),$ and let
$\overline{\Omega}$ denote the closure of the connected component of
$N-F,$ where $F$ is the fixed point set of $J,$ containing
$\infty_1.$
Then $X:\overline{\Omega}-\{\infty_1\}  \to  \l^3$ given by
$X(z):=\mbox{Real} \int^z \big(\phi_1,\phi_2,\phi_3 \big),$ where $\phi_1$
and $\phi_2$ are defined as asually, defines a  CMF graph with $n+1$
conelike singularities.
\end{corollary}

\begin{proof}
First, observe that the fixed point set of $J$ is not empty and consists of
the $n+1$ analytic circles $a_j:=z^{-1}(L_j),$ $j=0,\ldots,n,$ where
$L_j$ are the pairwise disjoint compact real intervals determined by the
points
$1,c_1,\ldots,c_n,-1,b_1,\ldots,b_n$ in $\r.$
By Koebe's Uniformization theorem, $\Omega$ is biholomorphic to a circular
domain with $n+1$ boundary components. 
It is straightforward to check that $(g,\phi_3)$ satisfies the hypothesis
of the Theorem \ref{th:exist}, and so $X$ defines a CMF graph with $n+1$
conelike singularities.
\end{proof}
The Riemann type graph in Figure \ref{fig:cat-rie} and the surfaces in Figure \ref{fig:2-3} lie in the family described in the previous corollary for $n=2$ and $n=3.$
\section{The moduli space of once punctured marked circular domains and its
associated bundles.} \label{sec:teich}
In  Section \ref{sec:moduli} we will prove the main theorems in the paper.
However, for a thorough explanation and subsequent development of these results  new tools and notations are required. This section is devoted to their introduction.

For any $c \in \c$ and $r>0,$ let $B_r(c)$ denote the round closed disc in
$\c$ of radius $r$
and centered at $c.$
In what follows, we will consider only circular domains bounded by $n+1$
circles, for a given $n\ge 1.$
Let $(\Omega-\{w\},a_0,\ldots,a_n)$ be an once punctured marked open circular domain, that
is to say, an open circular domain $\Omega\subset\overline{\c}$ punctured at $w \in \Omega$ together with an
ordering $(a_0,\ldots,a_n)$ of the circles in $\partial \Omega.$ 
Two once punctured marked open circular domains $(\Omega-\{w\},a_0,\ldots,a_n)$ and $(\Omega'-\{w'\},a'_0,\ldots,a'_n)$ are said to be equivalent if there exists a
biholomorphism $L:\Omega \to \Omega'$ (in fact, a Möbius transformation) such that $L(w)=w'$ and
$L(a_j)=a'_j,$ $j=0,\ldots,n.$
We call ${\cal T}_n$ the corresponding quotient space of 
equivalence classes of once punctured marked open circular domains, and refer to it as the  moduli space of once punctured marked circular domains with $n+1$ boundary components.

Given $[(\Omega-\{w\},a_0,\ldots,a_n)] \in {\cal T}_n,$ there exist unique
numbers
$c_1\in ]1,+\infty[,$ $c_2,\ldots,c_n\in\c$ and $r_1,\ldots,r_n\in\r^+$ and Möbius transformation $L$ in $\overline{\c}$ such that  $L(a_0)=\partial B_1(0),$
$L(a_j)=\partial B_{r_j}(c_j),$ $j=1,\ldots,n$ and $L(w)=\infty.$ Hence, $(\c-\Big( \cup_{j=0}^n B_{r_j}(c_j)\Big),\partial B_{1}(0), \partial B_{r_1}(c_1),\ldots,\partial B_{r_n}(c_n))$ is a representative of  $[(\Omega-\{w\},a_0,\ldots,a_n)],$ and  ${\cal T}_n$ can be canonically identified with the open subset in
$]1,+\infty[ \times \c^{n-1} \times (\r^+)^n \subset \r^{3n-1}$ consisting of those
points $v=(c_1,\ldots,c_n,r_1, \ldots,r_n)$ for which the balls $B_1(0)$
and $B_{r_j}(c_j),$ $j=1,\ldots,n,$ are pairwise disjoint.
\begin{figure}[htbp] 
\begin{center}
\includegraphics[width=12.6cm]{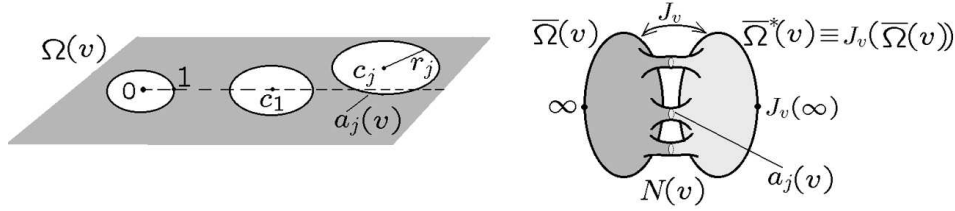}
\caption{$\Omega(v),$ $N(v)$ and $J_v.$ \label{fig:omega}}
\end{center}
\end{figure}

Given $v=(c_1,\ldots,c_{n},r_1, \ldots,r_{n}) \in {\cal T}_n,$ we call  
$\Omega(v)=\overline{\c}-\Big( \cup_{j=0}^n B_{r_j}(c_j)\Big),$  where $r_0=1,$ $c_0=0,$ and refer to the coordinates $c_i$ and $r_i$ of $v$  as $c_i(v)$ and
$r_{i}(v),$
respectively, $i=1,\ldots,n.$  Observe that the puncture corresponds to $\infty \in \Omega(v).$ We call $\overline{\Omega(v)}^*$ the mirror
(or double) of $\overline{\Omega(v)}=\Omega(v)\cup \big(\cup _{j=0}^n \partial B_{r_j(v)}(c_j(v))\big),$ with
the convention $c_0(v)=0$ and $r_0(v)=1,$ and put $N(v)=\overline{\Omega(v)} \cup \overline{\Omega(v)}^*$ for the associated closed Riemann surface.
Recall that $\overline{\Omega(v)}^* \cap \overline{\Omega(v)}$ consists of
the $n+1$ analytic circles $a_j(v):=\partial B_{r_j(v)}(c_j(v)),$ $j=0,\ldots,n.$
Moreover, we denote by $J_v:N(v)
\to N(v)$ the antiholomorphic involution  applying any point to its
mirror
image. Note that the fixed point set of $J_v$ coincides with $\cup_{j=0}^n
a_j(v).$

\begin{remark} \label{re:mirror}
A conformal model for $\overline{\Omega(v)}^*,$ $v \in {{\cal T}_n},$
consists of the planar domain $\overline{\Omega(v)}^*:=\{J_{v}(z) \;:\; z
\in \overline{\Omega(v)}\},$ where $J_{v}(z):=1/\bar{z}$ is the Schwarz
reflection about $\partial B_1(0).$
Moreover, $N(v)$ can be identified to the quotient of
$\overline{\Omega(v)} \cup \overline{\Omega(v)}^*$ under the
identification $z \equiv J_{v}(z),$ $z \in \partial \Omega(v).$

\end{remark}

\subsection{The bundles of divisors, meromorphic functions and meromorphic 1-forms}
Given $N$ a Riemann surface and $k$ a positive integer, we call
$$Div_{k}(N)=\{D \;:\; D \;\mbox{is an integral multiplicative divisor
on}\; N \; \mbox{of degree} \; k\}.$$ It well known that $Div_k(N)$ is the
quotient of $N^k$ under the action of the group of permutations of order
$k,$ and we denote by $p_k:N^k \to Div_k(N)$ the canonical projection. We
endow $Div_k(N)$ with the natural analytic structure induced by $p_k.$ More
precisely,
take $D_0=P_1^{m_1} \ldots P_s^{m_s} \in Div_k(N)$  and consider
$U=U_1^{m_1} \times \ldots \times U_s^{m_s} \subset N^k,$ where
$U_j^{m_j}=U_j \times \stackrel{m_j}{\ldots} \times U_j,$ $(U_j,z_j)$ is a
conformal chart around $P_j$ in $N,$ $z_j(P_j)=0,$ and $U_{j_1} \cap
U_{j_2}=\emptyset,$ $j_1 \neq j_2.$
The map $\xi:p_k(U) \to \c^k$ defined by $\xi (\prod_{j=1}^s Q_{1,m_j}
\cdot \ldots \cdot Q_{m_j,m_j})=((t_{1,m_j},\ldots,t_{m_j,m_j})_{j=1,
\ldots,s}),$ where $t_{h,m_j}=\sum_{l=1}^{m_j} (z_j(Q_{l,m_j}))^h,$
$h=1,\ldots,m_j,$ $j=1,\ldots,s$ defines an analytic parameterization
around $D_0.$ For more details, see \cite{farkas}.

The bundle of $(k_1,k_2)$-divisors over ${\cal T}_n$ is defined by $\div_{k_1,k_2}=\cup_{v
\in {{\cal T}_n}} \{(v,D_1,D_2) \;:\; D_i \in Div_{k_i}(\Omega(v)),\;
i=1,2\}.$  We endow $\div_{k_1,k_2}$ with its natural analytic structure.
To be more precise, let $(v_0,D_{1,0},D_{2,0}) \in \div_{k_1,k_2},$ take
$\epsilon>0$ small, and label $V(\epsilon)$ as the Euclidean ball of radius
$\epsilon$ in ${{\cal T}_n}$ centered  at $v_0.$
Write $D^i_0=z_{i,1}^{m_{i,1}} \ldots z_{i,s_i}^{m_{i,s_i}},$ $\sum_{h=1}^{s_i} m_{i,h}=k_i,$  $i=1,2.$
If $\epsilon$ is small enough, the set $W=\cap_{v \in V(\epsilon)}
\Omega(v)$ contains an open disc $U_{i,j}$ of radius $\epsilon$ around
$z_{i,j},$ $j=1,\ldots,s_i,$ $i=1,2,$ and
we can also take the discs $\{U_{i,j} \;:\; j=1,\ldots,s_i\}$
pairwise disjoint, $i=1,2.$
Consider the conformal charts $(U_{i,j},w_{i,j}:=z-z_{i,j}).$
Put $U_{i}=U_{i,1}^{m_{i,1}} \times \ldots \times U_{i,s_i}^{m_{i,s_i}},$
and observe that $U_i$ can be viewed as a subset of $\Omega(v)^{k_i},$ for
any $v
\in V(\epsilon).$ Likewise, $p_{k_i}(U_i) \subset Div_{k_i}(\Omega(v)),$
for any $v \in V(\epsilon),$ $i=1,2.$  The natural chart $\xi_i:p_{k_i}(U_{i,v}) \to
\c^{k_i}$ is
uniformly defined as before for any $v \in {{\cal T}_n},$ $i=1,2.$ By
definition ${\cal V}(\epsilon):=\{(v,D_1,D_2) \;:\; D_i \in
p_{k_i}(U_i)\subset
Div_{k_i}(\Omega(v)),\; i=1,2, \;v \in V(\epsilon)\}$ is a neighborhood of
$(v_0,D_{1,0},D_{2,0})$ in $\div_{k_1,k_2},$ and we say that $V(\epsilon)$
is its
associate open ball in ${{\cal T}_n}.$ Moreover, the map
$\Psi: {\cal V}(\epsilon) \to V(\epsilon) \times \xi_1(p_{k_1}(U_1)) \times
\xi_2(p_{k_2}(U_2))$ given  by:
\begin{equation} \label{eq:par1}
\Psi (v,D_1,D_2)=(v,\xi_1(D_1),\xi_2(D_2))
\end{equation}
defines a local (analytic)
parameterization around $(v_0,D_{1,0},D_{2,0})$ in $\div_{k_1,k_2}.$
We call $\vg:\div_{k_1,k_2} \to {{\cal T}_n},$ $\vg((v,(D_1,D_2)))=v,$
the natural projection.

In the sequel, and for the sake of simplicity, we simply write $\div_k$
instead of $\div_{k,0}$ or $\div_{0,k},$  and refer to
it as the bundle of k-divisors. We also establish the convention
$\div_{0,0}={\cal T}_n.$
\vspace{0.2cm}

For any $v \in {{\cal T}_n},$ call ${\cal C}(v)$ the family of meromorphic
functions on $N(v).$ The corresponding  bundle over ${\cal T}_n$ is denoted by ${\cal
C}_n=\cup_{v \in {{\cal T}_n}} {\cal C}(v).$

Likewise, we call  ${\cal H}(v)$  the space of  meromorphic $1$-forms on $N(v)$ and denote  by ${\cal H}_n=\cup_{v \in {{\cal
T}_n}} {\cal H}(v)$ the associated bundle over ${\cal T}_n.$

We need to introduce a convenient concept of  differentiability for maps
from $\div_{k_1,k_2}$ into ${\cal C}_n$ or
${\cal H}_n$ preserving the fibers. We start with the following definitions:
\begin{definition}
Let $M_j$ be a real manifold of dimension $m_j,$ $j=1,2,3,$ and let $f:M_1
\times M_2 \to M_3$ be a ${\cal C}^k$ map.
The map $f$ is said to be differentiable (or smooth) with $k$-regularity in $M_1$ if, for any
charts $(U_1 \times U_2,x\equiv
(x_1,\ldots,x_{m_1}),y\equiv(y_1,\ldots,y_{m_2}))$
in $M_1 \times M_2$ and $(U_3,z\equiv (z_1,  \ldots,z_{m_3}))$ in $M_3,$
the local expression of $f,$ $f(x,y):x(U_1)\times y(U_2) \to z(U_3),$
satisfies
that $f(\cdot,y)$ is smooth in $x(U_1)$ for any $y \in y(U_2),$ and all the
partial derivatives of $f(x,y)$ with respect to variables in $x$ are
${\cal C}^k$ in $x(U_1) \times y(U_2).$
\end{definition}
\begin{definition}
Let $v_0\in {\cal T}_n$ and $\epsilon >0$ small enough. Denote by
$V(\epsilon)$ the Euclidean ball of radius $\epsilon$ in ${\cal T}_n$
centered at $v_0.$
Since $V(\epsilon)$ is simply connected, standard homotopy arguments in
differential topology show the existence of a family of diffeomorphisms
$\{F_v:N(v_0) \to N(v)\;:\; v \in V(\epsilon)\}$ such that $F_{v_0}=\mbox{Id},$
$F_v(\infty)=\infty,$ $J_v\circ F_v\circ J_{v_0}=F_v,$ for any
$v \in V(\epsilon),$
and  $F:V(\epsilon)\times \overline{\Omega(v_0)}\to \c,$
$F(v,z):=F_v(z),$ is smooth.
By definition, we say that $\{F_v:N(v_0) \to N(v)\;:\; v \in V(\epsilon)\}$
is a smooth deformation of $N(v_0).$
Moreover note that, for $\epsilon$ small enough,  $\frac{\partial
F}{\partial z}\neq 0$ in $V(\epsilon)\times \overline{\Omega(v_0)}.$
\end{definition}

Let $W \subset \div_{k_1,k_2}$ be a submanifold, and let  $h: W \to {\cal
C}_n$ be a map preserving the fibers, that is to say,
$h_{v,D_1,D_2}:=h(v,D_1,D_2)\in{\cal C}(v)$ for any $(v,D_1,D_2)\in {W}.$
We are going to define the notion of differentiability with $k$-regularity of $h.$
Take ${\cal V}(\epsilon)$ any coordinate neighborhood in $\div_{k_1,k_2}$
defined as above and meeting ${W}.$
Denote by $V(\epsilon)$ the open ball in ${\cal T}_n$ associated to ${\cal
V}(\epsilon),$ and call $v_0  \in V(\epsilon)$ its center. Take a smooth
deformation of $N(v_0),$
$ \{F_v:N(v_0) \to N(v)\;:\; v \in V(\epsilon)\}.$
We say that $h$ is differentiable with $k$-regularity in ${\cal V}(\epsilon)\cap W$  if the map
$\hat{h}:({\cal V}(\epsilon) \cap W) \times N(v_0)\to \overline{\c},$ given by
$\hat{h}((v,D_1,D_2),x)=h_{v,D_1,D_2}(F_v(x))$ is smooth  with $k$-regularity in ${\cal
V}(\epsilon) \cap W$.  The map $h$ is said to be
differentiable with $k$-regularity on ${W}$  if it does in ${\cal V}(\epsilon) \cap
W,$ for any coordinate neighborhood ${\cal V}(\epsilon)$ meeting ${W}.$ It
is easy to check that this definition does not depend on choice of the
smooth deformation of $N(v_0).$

Let $\omega:{W}\to {\cal H}_n$ be a map preserving the fibers, that is  to
say, $\omega_{v,D_1,D_2}:=\omega(v,D_1,D_2)\in {\cal H}(v)$ for any
$(v,D_1,D_2)\in {W}.$
Take ${\cal V}_{\epsilon},$ $V(\epsilon),$ $v_0$ and $\{F_v:N(v_0) \to
N(v)\;:\; v \in V(\epsilon)\},$ as above and define
$\hat{\omega}:{\cal V}({\epsilon})\cap W\to {\cal H}(v_0)$ by
$\hat{\omega}(v,D_1,D_2)=\big( F_v^*(\omega_{v,D_1,D_2}) \big)^{(1,0)},$
where the superscript $(1,0)$ means the $(1,0)$ part of the 1-form
(by definition $(f\,dz+g\,d\overline{z})^{(1,0)} =f\,dz$).
We say that $\omega$ is differentiable in with $k$-regularity ${\cal V}(\epsilon)\cap W$ 
if for any local chart $(U,z)$ in $N(v_0),$ the map $\hat{f}:({\cal
V}({\epsilon})\cap W)\times U\to \overline{\c},$ given by
$\hat{f}((v,D_1,D_2),z)=\hat{\omega}(v,D_1,D_2)(z)/dz$ is smooth with $k$-regularity in ${\cal
V}(\epsilon)\cap W.$ The global concept of
differentiability with $k$-regularity in ${W}$ is defined in the obvious way.\\

\subsection{The Jacobian bundle}
In order to define the Jacobian bundle over ${\cal T}_n,$ some topological and analytic
preliminaries are required.

Take $v\in{\cal T}_n,$ and whenever no confusion is
possible, identify the homology classes of the boundary circles $a_j(v),$ $j=0,1,\ldots,n,$ with their representing curves.

Let $b_1,\ldots,b_n$ be closed curves in $N(v)$ such that
$B=\{a_1(v),\ldots,a_n(v),b_1,\ldots,b_n\}$ is a canonical homology basis,
that is to say, the intersection numbers $(a_j(v),a_h(v)),$ $(b_j,b_h)$
vanish, and $(a_j(v),b_h)=\delta_{ah},$ where $\delta_{ah}$ refers to the
Kronecker symbol.

It is possible to make a canonical choice of the homology basis of $N(v).$
Indeed, note first that $J_v$ fixes $a_j(v)$ pointwise, and so,
$J_v(a_j(v))=a_j(v).$  Take a curve $\gamma_j \subset \overline{\Omega(v)}$
joining $a_0(v)$ to $a_j(v),$ and observe that the curve $b_j(v)$ obtained
by joining $\gamma_j$ and $J_v(\gamma_j)$  satisfies $J_v(b_j(v))=-b_j(v)$
in the homological sense, and its homology class does not
depend on the choice of $\gamma_j.$ In other words, the identity $J_v(b_j(v))=-b_j(v)$
characterizes
$B(v)=\{a_1(v),\ldots, a_n(v),b_1(v),\ldots,b_n(v)\}$ as canonical homology
basis of $N(v).$

Let $\{\eta_1(v),\ldots,\eta_n(v)\}$ be the dual basis of $B(v)$ for the
space of holomorphic 1-forms on $N(v),$
that is to say, the unique basis satisfying
$\int_{a_k(v)}\eta_j(v)=\delta_{jk},$ $j,k=1,\ldots,n.$
Call $\Pi(v)=(\pi_{j,k}(v))_{j,k=1,\ldots,n}$ the associated matrix of
periods, $\pi_{j,k}(v)=\int_{b_j(v)} \eta_k(v).$

Given $D=\prod_{j=1}^s w_j^{m_j} \in Div_k(\Omega(v)),$ we
denote by $\tau_D(v)$ the unique meromorphic 1-form on $N(v)$ having simple
poles at $w_j$ and $J_v(w_j),$ $j=1,\ldots,s,$ and no other poles,
and satisfying $\mbox{Residue}_{w_j} \big(\tau_D(v)\big)=
-\mbox{Residue}_{J_v(w_j)}\big(\tau_D(v)\big)=-m_j,$ $\int_{a_i(v)}
\tau_D(v)=0,$ for any $j,i.$
Likewise, take $D_1 =\prod_{j=1}^s w_{j,1}^{m_j},D_2=\prod_{h=1}^r
w_{h,2}^{n_h} \in Div_k(\Omega(v))$ and define  $\kappa_{D_1,D_2}(v)$ as
the unique
meromorphic 1-form on $N(v)$ having simple poles at $w_{j,1},\,w_{h,2}$ and
$J_v(w_{j,1}),\, J_v(w_{h,2}),$ $j=1,\ldots,s,$ $h=1,\ldots,r,$ and no other poles, and
satisfying
$$\mbox{Residue}_{w_{j,1}}\big(\kappa_{D_1,D_2}(v)\big)=
\mbox{Residue}_{J_v(w_{j,1})} \big(\kappa_{D_1,D_2}(v)\big)=-m_j,$$ $$ 
\mbox{Residue}_{w_{h,2}}\big(\kappa_{D_1,D_2}(v)\big)=
\mbox{Residue}_{J_v(w_{h,2})}\big(\kappa_{D_1,D_2}(v)\big)=n_h$$
and $\int_{a_i(v)} \kappa_{D_1,D_2}(v)=0,$ for any $j,h,i.$

Our aim is to show that $\eta_j(v),$ $\tau_D(v)$ and $\kappa_{D_1,D_2}(v)$
depend smoothly with $1$-regularity on $v,$ $(v,D)$ and $(v,D_1,D_2),$ respectively. This fact
is enclosed in the following technical lemma.

\begin{lemma} \label{lem:schauder}
Let $\Omega$ be an open bounded domain in $\r^2$ with smooth boundary and
$\b$ an open Euclidean ball in $\r^m.$ Let
$H(t,x),\phi(t,x):\overline{\b} \times \overline{\Omega} \to \r$ be two ${\cal
C}^1$ functions, where
$\phi_t:=\phi(t,\cdot) \in {\cal C}^{2,\alpha}(\overline{\Omega}),$ $
\alpha \in ]0,1],$ for all $t=(t_1,\ldots,t_k),$ and the map $t \mapsto
\phi_t$ is ${\cal
C}^1$ in the ${\cal C}^{2,\alpha}$ norm on $\overline{\Omega}.$ Consider a
smooth one parameter family of metrics $ds^2_{t}$ on $\overline{\Omega},$
$t\in
\overline{\b}$ and denote by $\Delta_t$ the associated  family of
Laplacians. Let $u_t \in {\cal C}^0(\overline{\Omega}) \cap {\cal
C}^2(\Omega)$ be the
solution of the boundary value problem $\Delta_{t} u_{t}= H_{t},$
$u_{t}|_{\partial {\Omega}} = \phi_{t}|_{\partial{\Omega}},$
where $H_{t}(x):=H(t,x),$  $(t,x) \in \overline{\b} \times
\overline{\Omega} .$

If we define $u:\b\times \overline{\Omega}\to \r$ by $u(t,x)=u_{t}(x),$ then:
\begin{enumerate}
\item $u_t \in {\cal C}^{2,\alpha}(\overline{\Omega}),$ and
\item the map $t \mapsto u_t$ is ${\cal C}^1$  in the ${\cal C}^{2,\alpha}$
norm on $\overline{\Omega},$
\item if $H$ and $\phi$ are smooth, then for any $q, m_1,\ldots,m_p\in \n,$
$m_1+\ldots+m_p=q,$ $i_1,\ldots,i_p \in \{1,\ldots,m\},$ the function
$\frac{\partial^q u_t }{\partial t_{i_1}^{m_1}\ldots\partial t_{i_p}^{m_p}}
\in {\cal C}^{2,\alpha}(\overline{\Omega}),$
\item if $H$ and $\phi$ are smooth, then the map $t \mapsto u_t$ is ${\cal
C}^\infty$   in the ${\cal C}^{2,\alpha}$ norm on $\overline{\Omega}.$ As a
consequence, for any $q, m_1,\ldots,m_p\in \n,$ $m_1+\ldots+m_p=q,$
$i_1,\ldots,i_p \in \{1,\ldots,m\},$
the function $\frac{\partial^q u }{\partial t_{i_1}^{m_1}\ldots\partial
t_{i_p}^{m_p}} \in {\cal C}^{2,\alpha}( (\b \times \Omega) \cup T),$ where
$T$ is any differentiable portion of $\partial (\b \times \Omega).$
\end{enumerate}
\end{lemma}

\begin{proof}
It is enough to consider the case $m=2$ (the general case is similar).

First, note that $(1)$ is a straightforward consequence of global
regularity theorem \cite{gilbarg} p.106.

The maximum principle and the classical Schauder estimates (\cite{gilbarg}
p.35 and p. 93), show that  the family $\{u_{t,s} \;:\; (t,s) \in
\overline{\b}\}$ is
bounded in the ${\cal C}^{2,\alpha}$ norm on $\overline{\Omega}.$  Fix
$(t_0,s_0),$ then for each $(t,s),$ the function $u_{t,s} -u_{t_0,s_0}$
satisfies:
$\Delta_{t_0,s_0}(u_{t,s}-u_{t_0,s_0})=
(\Delta_{t_0,s_0}-\Delta_{t,s})u_{t,s} +H_{t,s} -H_{t_0,s_0}$ on $\Omega$
and $u_{t,s}-u_{t_0,s_0} = \phi_{t,s} -
\phi_{t_0,s_0}$ on $\partial\Omega.$ The maximum principle and Schauder's
estimates then show that the map $(t,s) \mapsto u_{t,s}$ is continuous at
$(t_0,s_0)$ with respect to the  ${\cal C}^{2,\alpha}$ norm on
$\overline{\Omega}.$  Let us show that the functions $w_{t,s_0}= (u_{t,s_0}
-u_{t_0,s_0})/(t-t_0)$ converge in the $C^{2,\alpha}$ norm on
$\overline{\Omega},$ as $t\to t_0,$ to the solution $y_{t_0,s_0}$ of the
problem:
$\Delta_{t_0,s_0} y_{t_0,s_0}= - \frac{\partial\Delta_t}{\partial
t}(t_0)(u_{t_0,s_0})+ \frac {\partial H_{t,s}}{\partial t} (t_0,s_0)$ on
${\Omega}$ and
$y_{t_0,s_0} = \frac{\partial \phi}{\partial t}(t_0,s_0)$ on $\partial
{\Omega}.$ Indeed, put
$L(t,s_0)= (\Delta_{t,s_0} - \Delta_{t_0,s_0})/(t-t_0)$ for $t\neq t_0.$
Then the functions $w_{t,s_0}$ are solutions
of the problem: $\Delta_{t_0,s_0}w_{t,s_0} = - L(t,s_0)(u_{t,s_0}) + \frac
{H_{t,s_0} - H_{t_0,s_0}}{ t-t_0}$ on ${\Omega},$
$w_{t,s_0}|_{\partial {\Omega}} =(\phi_{t,s_0} -\phi_{t_0,s_0})/(t-t_0).$
On the other hand, the function $w_{t,s_0} - y_{t_0,s_0}$ satisfies:
$\Delta_{t_0,s_0}(w_{t,s_0} - y_{t_0,s_0}) = - L(t,s_0)u_{t,s_0} +
\frac{\partial\Delta_{t,s}}{\partial t}(t_0,s_0)u_{t_0,s_0}+ \frac
{H_{t,s_0} - H_{t_0,s_0}}{( t-t_0)}
-\frac {\partial H_{t,s}}{\partial t} (t_0,s_0),$
$(w_{t,s_0} - y_{t_0,s_0})|_{\partial {\Omega}}=(\phi_{t,s_0}
-\phi_{t_0,s_0})/(t-t_0)-\frac{\partial \phi}{\partial t}(t_0,s_0).$
Therefore, as before, the maximum principle and Schauder's estimates imply
that $w_{t,s_0}$ converges  to $y_{t_0,s_0}$ in the $C^{2,\alpha}$ norm on
$\overline{\Omega}.$
Likewise,  the family $\{y_{t,s} \;:\; (t,s) \in \overline{\b}\}$ is
bounded with respect to the ${\cal C}^{2,\alpha}$ norm on
$\overline{\Omega}.$ Furthermore, the function $y_{t_1,s_1}-y_{t_2,s_2}$
satisfies the equation $\Delta_{t_1,s_1} (y_{t_1,s_1} - y_{t_2,s_2}) = -
\frac{\partial{\Delta }}{\partial t}(t_1,s_1)u_{t_1,s_1}
+\frac{\partial{\Delta }}{\partial t}(t_2,s_2)u_{t_2,s_2} + \frac{\partial
H_{t,s}}{\partial t} (t_1,s_1) - \frac{\partial H_{t,s}}{\partial t}
(t_2,s_2) - (\Delta_{t_1,s_1} -
\Delta_{t_2,s_2}) y_{t_2,s_2},$ $(y_{t_1,s_1} - y_{t_2,s_2})|_{\partial
{\Omega}}=\frac{\partial \phi}{\partial t}(t_1,s_1)-\frac{\partial
\phi}{\partial t}(t_2,s_2),$ and hence, using again  the maximum principle
and Schauder's estimates, we obtain the continuity
of $y_{t,s}$ in $(t,s)$ in the $C^{2,\alpha}$ norm on $\overline{\Omega}.$
The same holds for $\frac{\partial u}{\partial s},$ and proves (2).

Suppose now that $H$ and $\phi$ are ${\cal C}^\infty.$ The above argument
can be applied
to $\frac{\partial u_{t,s}}{\partial t}$ and $\frac{\partial
u_{t,s}}{\partial s},$ and so $(1)$ and $(2)$ also hold for these
functions. An iterative argument proves that the map $(t,s) \mapsto
u_{t,s}$ is ${\cal C}^\infty$   in the ${\cal C}^{2,\alpha}$ norm on
$\overline{\Omega},$ which proves $(3)$ and the first part of $(4).$ For
the second part of $(4),$ let $f(t,s,x)$ denote any partial derivative of
$\frac{\partial^q u }{\partial t^{m_1} \partial s^{m_2}},$ $m_1+m_2=q,$ $q \in \n,$ of
order two with respect to variables in $\b \times \overline{\Omega}.$ It is enough to
check that $||f||_{0,\alpha}$ (where $||\cdot||_{0,\alpha}$ is the ${\cal
C}^{0,\alpha}$ norm in $\overline{\b} \times \overline{\Omega}$) is
bounded. This follows from the inequality:
$$||f||_{0,\alpha} \leq C \left(\mbox{Max} \{||f(t,s,\cdot)||_{0,\alpha}
\;:\; (t,s) \in \b\}+ ||\frac{\partial f}{\partial t}||_0+ ||\frac{\partial
f}{\partial s}||_0  \right),$$ where $C$ is a positive constant and
$||\cdot||_0$ is the ${\cal C}^0$ norm on $\overline{\b} \times
\overline{\Omega}.$ The ${\cal C}^{2,\alpha}$ regularity of $u$ and all its
partial derivatives in $(t,s)$ on the smooth portion of $\partial(\b \times
\Omega)$ follows also from the regularity theorem.

\end{proof}

\begin{corollary} \label{co:smooth}
The maps $\eta_j:{{\cal T}_n}\to {\cal H}_n,$ $v\mapsto\eta_j(v),$
$\tau:\div_{k}\to{\cal H}_n,$ $(v,D)\mapsto\tau_D(v),$ and
$\kappa:\div_{k,k}\to{\cal H}_n,$ $(v,D_1,D_2)\mapsto\kappa_{D_1,D_2}(v)$
are differentiable with $1$-regularity.

As a consequence, the functions $\pi_{j,k}(v):= \int_{b_j(v)} \eta_k(v),$
are differentiable on ${{\cal T}_n}.$
\end{corollary}

\begin{proof}
Let $v_0 \in {{\cal T}_n},$ $(v_0,D_0) \in \div_k$ and
$(v_0,D_{1,0},D_{2,0}) \in \div_{k,k},$
and take $V(\epsilon),$ ${\cal V}_0(\epsilon)$ and ${\cal V}(\epsilon)$
the  previously defined open neighborhoods of $v_0, $ $(v_0,D_0)$ and
$(v_0,D_{1,0},D_{2,0})$ in ${{\cal T}_n},$ $\div_k$ and $\div_{k,k},$
respectively.

Write $D =\prod_{j=1}^l w_{j}^{m_j},$ $D_1 =\prod_{j=1}^{l_1}
w_{j,1}^{m_{j,1}}$ and $D_2=\prod_{h=1}^{l_2} w_{h,2}^{m_{h,2}},$ and
denote by $A_{v,D}= \sum_{j=1}^l \log |{z-w_j}|^{m_j},$ $A_{v,D_1,D_2}=
\mbox {Im} \left(\int \big( \sum_{h=1}^{l_2} \frac{m_{h,2}}{z-w_{h,2}}-
\sum_{j=1}^{l_1} \frac{m_{j,1}}{z-w_{j,1}} \big) dz \right).$ Observe that
$A_{v,D}$ is well defined on
$\overline{\Omega(v)}-\{\infty,w_{1},\ldots,w_{l}\},$  and $A_{v,D_1,D_2}$
is well defined in a small enough neighborhood of $\partial \Omega(v)$
consisting of the union of $n+1$ small annuli, up to adding constants.

Let $h_{j,v}$ be the unique harmonic function on
$\Omega(v)$ satisfying $h_{j,v}|_{a_k(v)}=\delta_{jk}.$
Call  also $h_{v,D}$ (resp. $h_{v,D_1,D_2}$) the unique
harmonic function on $\Omega(v)$ such that
$h_{v,D}|_{\partial \Omega(v)}=A_{v,D}$
(resp. $h_{v,D_1,D_2}|_{\partial \Omega(v)}=A_{v,D_1,D_2}$).

Let $\hat{\eta}_j(v),$ $\hat{\tau}_D(v),$ $\hat{\kappa}_{D_1,D_2}(v)$
denote the 1-forms $\partial_z h_{j,v},$ $2 \partial_z \big(
h_{v,D}-A_{v,D} \big)$ and  $2 i \partial_z \big(h_{v,D_1,D_2}-A_{v,D_1,D_2}
\big),$ which are well defined as meromorphic 1-forms on $\Omega(v).$ They are extended by Schwarz reflection and with the same name to
$N(v).$ Moreover,  $\{\hat{\eta}_j(v) \;:\;j=1,\ldots,n\}$ is a basis
of the complex linear space of holomorphic 1-forms on $N(v),$ and
$\hat{\tau}_D(v),$ $\hat{\kappa}_{D_1,D_2}(v)$ are meromorphic 1-forms
having the same poles (with the same residues) as ${\tau}_D(v)$ and
${\kappa}_{D_1,D_2}(v),$ respectively.

\vspace{0.3cm}

\begin{quote} {\bf Claim:} {\em The maps $\hat{\eta}_j:V(\epsilon)\to {\cal
H}_n,$ $v\mapsto\hat{\eta}_j(v),$
$\hat{\tau}:{\cal V}_0(\epsilon)\to{\cal H}_n,$
$(v,D)\mapsto\hat{\tau}_D(v),$ and
$\hat{\kappa}:{\cal V}(\epsilon)\to{\cal H}_n,$
$(v,D_1,D_2)\mapsto\hat{\kappa}_{D_1,D_2}(v),$ are smooth with $1$-regularity.}

\end{quote}

Take a  smooth deformation $\{F_v:N({v_0}) \to N({v}), \; v \in
V(\epsilon)\}$ of $N(v_0).$ Note that $F_v(a_j(v_0))=a_j(v),$
$F_v(b_j(v_0))=b_j(v),$ in the homological sense, for any $v \in
V(\epsilon)$ and any $j.$

In the sequel, and for the sake of simplicity, we will put $h_{j,v}=h_v.$

Let $\Gamma:\overline{\b} \to \overline{W}$ be a parameterization in ${\cal
T}_n,$ $\div_k$ or
$\div_{k,k}$, where $W$ is an open neighborhood contained in either $V(\epsilon),$
${\cal V}_0(\epsilon)$ or ${\cal V}(\epsilon),$  and write
$\Gamma(t)=v(t),$ $\Gamma(t)=(v(t),D(t))$ or
$\Gamma(t)=(v(t),D_1(t),D_2(t))$ in each case, $t \in \b.$
Call $F_t:=F_{v(t)}:\overline{\Omega(v_0)}\to\overline{\Omega(v(t))}$ and
$h_t:=h_{\Gamma(t)}:\overline{\Omega(v(t))}\to\r.$ Then, it suffices to
check that the map
$\hat{u}:\overline{\b}\times\overline{\Omega(v_0)}\to \r,$ $(t,x)\mapsto
h_t(F_t(x)),$ is smooth with $2$-regularity in $\b.$ 

To do this, observe first that 
$\Omega(v_0)\subset\overline{\c}$ is conformally equivalent to the bounded
domain $\Omega'=\{1/x\;:\;x\in\Omega(v_0)\} \subset \r^2,$ where the biholomorphism is given by $T:\overline{\Omega'}\to\overline{\Omega(v_0)},$
$T(x)=1/x.$ If we put
$u:\overline{\b}\times\overline{\Omega'}\to\r,$ $u(t,x)=\hat{u}(t,T(x)),$ it is clear that $\hat{u}$ is smooth with $2$-regularity in $\b$  if and only if
$u$ does. 

Consider now the metric $ds_t^2$ on
$\overline{\Omega'}$ making
$F_t\circ T:\overline{\Omega'}\to \overline{\Omega(v(t))}$ an isometry, and
denote by
$\Delta_t$ the associated family of Laplacians, $t\in\overline{\b}.$
Then $u_t:=h_t(F_t(T(x)))$ is the solution of the boundary value problem
$\Delta_t u_t=0,$ $u_t|_{\partial\Omega'}=\phi_t|_{\partial\Omega'},$
$t\in\overline{\b},$  where $\phi_t(x)=\phi(t,x)$ and $\phi$ is a suitable smooth function in
$\overline{\b}\times\overline{\Omega'}.$ From Lemma \ref{lem:schauder}, $u$
and its partial derivatives till the second order are smooth with $2$-regularity in $\b,$ which proves the claim.

\vspace{0.3cm}

From the previous claim, we infer that the maps $\hat{\eta}_j:{{\cal
T}_n}\to {\cal H}_n,$ $v\mapsto\hat{\eta}_j(v),$
$\hat{\tau}:\div_{k}\to{\cal H}_n,$ $(v,D)\mapsto\hat{\tau}_D(v),$ and
$\hat{\kappa}:\div_{k,k}\to{\cal H}_n,$
$(v,D_1,D_2)\mapsto\hat{\kappa}_{D_1,D_2}(v)$ are smooth with $1$-regularity.
Hence the period functions $v \mapsto \int_{a_h(v)} \hat{\eta}_j(v)$ are
smooth on ${{\cal T}_n},$ and since $\hat{\eta}_j(v)= \sum_{h=1}^n \big(
\int_{a_h(v)} \hat{\eta}_j(v)\big) \eta_h(v),$ for all $j,$ we easily check
that $\eta_1,\ldots, \eta_n$ are smooth with $1$-regularity in ${{\cal T}_n}$ .

Moreover, $\tau_D(v)=\hat{\tau}_D(v)-  \sum_{h=1}^n \big( \int_{a_h(v)}
\hat{\tau}_D(v)\big) \eta_h(v),$  and likewise
$\kappa_{D_1,D_2}(v)=\hat{\kappa}_{D_1,D_2}(v)-  \sum_{h=1}^n \big(
\int_{a_h(v)} \hat{\kappa}_{D_1,D_2}(v)\big) \eta_h(v).$  Reasoning as
before, $\tau$
and $\kappa$ are smooth with $1$-regularity in $\div_k$ and $\div_{k,k},$ respectively, which
concludes the proof.

\end{proof}

For any $v \in {{\cal T}_n},$ let $L(v)$ be the lattice over $\z$ generated
by $\{e^1,\ldots,e^n,\pi^1(v),\ldots,\pi^n(v)\},$ where $e^j=\,^T(0,\ldots,
\stackrel{j}{1},\ldots,0)$ and $\pi^j(v)=\,^T(\pi_{1,j}(v),
\ldots,\pi_{n,j}(v)).$ Obviously, $L(v)$ depends smoothly on $v.$
Put ${\cal J}(v)=\c^n/L(v)$  for the Jacobian variety associated to $N(v)$
and label $\pg_v:\c^n \to {\cal J}(v)$ as the natural projection.

Denote by ${{\cal J}}_n=\cup_{v \in {{\cal T}_n}} \{(v,q) \;:\; q \in
{\cal J}(v)\}$ the Jacobian bundle and define $\pg:{{\cal T}_n} \times \c^n
\to {{\cal J}}_n$ by $\pg(v,z)=(v,\pg_v(z)).$ Since $\pg$ is locally
injective, ${{\cal J}}_n$  can be endowed with the analytic structure
making $\pg$ a local diffeomorphism.
To be more precise, let $(v_0,q_0) \in {{\cal J}}_n,$ and consider $w_0 \in
\c^n$ such that $\pg_{v_0}(w_0)=q_0.$ Let $W(\epsilon)$ be an open ball of
radius $\epsilon$ in $\c^n$ centered at $w_0$ such that
$\pg_{v_0}|_{W(\epsilon)}:W(\epsilon)  \to \pg_{v_0}(W(\epsilon))$ is a
conformal diffeomorphism.
Since $L(v)$ depends smoothly on $v,$ $\pg_{v}|_{W(\epsilon)}:W(\epsilon)
\to \pg_{v}(W(\epsilon))$ is also a conformal diffeomorphism for any $v$ in
the open ball $V(\epsilon)$ of radius $\epsilon$ centered at $v_0$ in
${{\cal T}_n},$  provided that $\epsilon$ is small enough.
By definition, the set
${\cal W}(\epsilon):=\pg(U(\epsilon)\times W(\epsilon))$ is a neighborhood
of $(v_0,q_0)$ in ${{\cal J}}_n,$ and the map
\begin{equation} \label{eq:par2}
\Upsilon:U(\epsilon) \times W(\epsilon) \to {\cal W}(\epsilon), \quad
\Upsilon=\pg|_{U(\epsilon) \times W(\epsilon)}
\end{equation}
is a local analytic parameterization.

\subsection{The n-spinorial bundle} \label{subsec:spin}

In this subsection we introduce a  bundle ${\cal S}_n$ over ${\cal T}_n$ of spinorial type. It  will arise in a natural way later when we study the space of CMF graphs with $n+1$ singular points.
From the technical point of view, ${\cal S}_n$ is a submanifold of
the bundle $\div_n$ of integral divisors of degree $n$ over ${\cal T}_n.$

Some terminology about the Jacobian bundle and its canonical and spinorial
sections is first required.

Consider the holomorphic 1-form $\overline{J_v^*(\eta_j(v))}.$ Taking into
account that $J_v$ fixes $a_j(v)$ pointwise, we infer that $\int_{a_k(v)}
\overline{J_v^*(\eta_j(v))} =\delta_{jk},$ and so,
$J_v^*(\eta_j(v))=\overline{\eta_j(v)}.$  Moreover, since
$J_v(b_j(v))=-b_j(v),$ then $\pi_{j,k}(v)=\int_{b_k} \eta_j(v)$ is an
imaginary number, for any $j$ and $k.$ It follows that there exists a
unique mirror involution $I_v:J(v) \to J(v)$ satisfying $I_v (\pg_v(w))=
\pg_v (\overline{w}),$ for any $w \in \c^n.$
This allows us to define the analytic involution ${\cal I}:{{\cal J}}_n
\to {{\cal J}}_n,$ ${\cal I}((v,q))=(v,I_v(q)),$ referred to as the mirror
involution in the Jacobian bundle.

For any $v \in {{\cal T}_n},$ we call $\varphi_v:N(v) \to {\cal J}(v)$ the
Abel-Jacobi embedding defined by
$$\varphi_v(z)=\pg_v\left( \int_{1}^{z} {^T(\eta_1(v),\ldots,\eta_n(v))}
\right),$$
Recall that $1\in \overline{\Omega(v)} \subset N(v)$ uniformly on $v,$ and since $J_v(1)=1,$
then $\varphi_v \circ J_v= I_v \circ \varphi_v,$ $v\in {\cal T}_n.$
We extend $\varphi_v$ with the same name to the Abel-Jacobi map
$\varphi_v:Div_k(N(v)) \to {\cal J}(v)$ given by $\varphi_v(P_1 \cdot
\ldots \cdot P_k)= \sum_{j=1}^k \varphi_v(P_j),$ $k \geq 1.$
We also denote by $\varphi:\div_{k} \to
{\cal J}_n,$ the map $\varphi(v,D)=(v,\varphi_v(D)).$
It is clear from Corollary \ref{co:smooth} that $\varphi$  is smooth.

Let $v \in {{\cal T}_n}$ and let $\omega$ be a non-zero meromorphic 1-form on
$N(v).$ Denote by $(\omega)$  the canonical divisor  associated to
$\omega,$  and call
$T(v)= \varphi_v( (\omega) )\in J(v) ,$ which by Abel's theorem does not depend on the
choice of $\omega.$ We call $\hat{T}:{\cal T}_n\to
\Jg_n,$ $\hat{T}(v)=(v,T(v)),$ the corresponding section of ${\cal J}_n.$
Taking for instance $\omega=\eta_j(v),$ which satisfies $J_v^*(\eta_j(v))=\overline{\eta_j(v)},$ $v \in {\cal T}_n,$  we can check that the divisor
$(w)$ is invariant under $J_v,$ and hence ${\cal I}(\hat{T}(v))=\hat{T}(v).$
A remarkable fact is that $\hat{T}$ is smooth. Indeed,  $T(v)=-2 K(v),$ where $K(v)$ is the vector of Riemann constants in $N(v)$ (see \cite{farkas} p. 298 (1980)). It is well known (see \cite{farkas} p. 290, 1980) that $$K(v)=\pg_v\left(\sum_{j=1}^n \big( (\frac{\pi_{j,j}(v)}{2})e^j-\int_{a_j(v)} \tilde{\varphi}_v \eta_j(v)\big)\right),$$ where $\tilde{\varphi}_v|_{a_j(v)}$ is any  lift to $\c^n$ with respect to $\pg_v$ of $\varphi_v|_{a_j(v)}.$ Since the map $\hat{K}:V_n \to {{\cal J}}_n,$ $v \mapsto (v,K(v)),$ is smooth, the same holds for $\hat{T}.$

A different approach to the regularity of $\hat{T}$ can be found after the proof of  Lemma \ref{lem:porfin}.

As a  consequence of the smoothness of $\hat{T},$ there are exactly $2^{2n}$
differentiable maps $\hat{K}_1,\ldots, \hat{K}_{2^{2n}}:{{\cal T}_n} \to {{\cal
J}}_n,$ $\hat{K}_j(v)=(v,K_j(v)),$ satisfying $2 K_j=T,$ for any $j.$ The following lemma shows that these {\em spinor sections} of the Jacobian bundle are invariant under the
mirror involution. 

\begin{lemma} \label{lem:spinor}
${\cal I} \circ \hat{K}_j=\hat{K}_j,$ for any $j=1,\ldots,2^{2 n}.$
\end{lemma}
\begin{proof}
Indeed, note that $I_v(K_j(v))=K_j(v)+\pg_v(\frac{1}{2} \sum_{h=1}^n
(m_h(v) e^h+n_h(v) \pi^h(v))),$ where $m_h(v),$ $n_h(v) \in \z$ are
continuous functions of $v.$ Using that ${{\cal T}_n}$ is connected we get
that $m_h(v),\;n_h(v)$ are constant.
Hence, the set ${\cal A}_j:= \{v \in {{\cal
T}_n} \;:\; I_v(K_j(v))=K_j(v)\}$ is either empty or the whole of ${{\cal
T}_n}.$ On the other hand, $K_j(v)=K_1(v)+q_j(v),$ where $2 q_j(v)=0,$ and so,
$I_v(q_j(v))=q_j(v).$  Therefore ${\cal A}_1={{\cal T}_n}$ if and only if
${\cal A}_j={{\cal T}_n}$ for any $j.$ Finally, consider the compact genus
$n$
Riemann surface $N=\{(z,w) \in \overline{\c}\;:\; w^2=\prod_{i=1}^{2n+2}
(z-c_i)\},$ where $c_i \in \r$ and $c_1<c_2< \ldots <c_{2n+2}.$ Then define
the
antiholomorphic involution $J(z,w)=(\overline{z},-\overline{w})$ and the
holomorphic 1-form $\omega=\prod_{i=1}^{n-1} (z-c_i) \frac{dz}{w}.$  The
function
$w$ has a well defined branch $w_+$ on the planar domain $\Sigma
=\overline{\c}-
\cup_{i=0}^{n} [c_{2 i+1},c_{2i+2}],$ and  the domain
$\{(z,w_+(z)) \;:\; z \in \Sigma\} \subset N$ is biholomorphic to a circular
domain $\Omega(v_0),$ $v_0 \in {{\cal T}_n}.$ Furthermore, up to this
biholomorphism, $N=N(v_0)$ and $J=J_{v_0}.$

Observe that the canonical divisor $(\omega)$ is given by $c_1^2 \cdot
\ldots \cdot c_{n-1}^2,$ where, up to the above identification $(c_i,0)
\equiv c_i \in N(v_0).$  Since $J_{v_0}(c_i)=c_i,$ then
$k_0:=\sum_{i=1}^{n-1} \varphi_{v_0} (c_i) \in J({v_0})$ is invariant under
$I_{v_0}$ and  $2 k_0=T(v_0).$ Up to relabeling, we can suppose
that $k_0=K_1(v_0)$ and hence ${\cal A}_1={{\cal T}_n}.$ This completes the
proof.
\end{proof}

The n-spinorial bundles are defined by
${\cal S}_n(i)=\{(v,D) \in \div_n \;:\;\varphi_v(D)-\varphi_v(J_v(\infty))=\hat{K}_i(v)\}.$  
We also call
$${\cal S}_n=\cup_{i=1}^{2^{2 n}} {\cal S}_n(i)=
\{(v,D)\in\div_n \;:\; \varphi_v(D^2)-\varphi_v(J_v(\infty)^2)=\hat{T}(v)\}, $$
and refer to it as  the (global) n-spinorial bundle.

\begin{theorem}[Structure] \label{th:submanifold}
The space ${\cal S}_n$ is a smooth real $(3n-1)$-dimensional submanifold of
$\div_{n}$
and the map $\vg:{\cal S}_n \to {\cal T}_n,$ $\vg(v,D) = v,$ is
a finite covering. 
\end{theorem}
\begin{proof}
The fact  ${\cal S}_n \neq \emptyset$ is a consequence of Corollary \ref{co:exist} and Proposition \ref{pro:cale}. The key step of this proof is that ${\cal S}_n$ does not contain any {\em special} divisor (see \cite{farkas}).

Consider the differentiable map $H: \div_{n} \to {\cal J}_n$ given by
$H(v,D)=(v,2\varphi_v(D)-2\varphi_v(J_v(\infty))-T(v)),$ and note that
${\cal S}_n=\{(v,D) \in \div_{n} \;:\;H(v,D)=(v,0)\}.$
In order to prove that ${\cal S}_n$ is a differentiable submanifold of
$\div_{n},$
it suffices to check that $dH_q$ is bijective at any point $q$ of ${\cal
S}_n.$

Let $q_0:=(v_0,D_0)$ be an arbitrary point of ${\cal S}_n.$
We are going to write the expression of $H$ in local coordinates $\Psi$
around $q_0 \in \div_{n}$ and $\Upsilon^{-1}$ around $H(q_0)=(v_0,0) \in
{\cal J}_n$ (see equations (\ref{eq:par1}) and (\ref{eq:par2})).

To do this, take $V(\epsilon),$ $V(\epsilon')$ balls in ${\cal T}_n$
centered at $v_0$ of radius $\epsilon$ and $\epsilon'$ respectively and
$W(\epsilon') \in
\c^n$ the ball of radius $\epsilon'$ centered at the origin.  Write
$D_0=z_1^{m_1} \cdot \ldots \cdot z_s^{m_s} \in Div_n(\Omega(v_0)),$
and denote by $(U_j,w_j:=z-z_j)$ the conformal parameter in $\Omega(v_0),$
where $U_j$ is the open disc of radius $\epsilon$ centered at $P_j,$
$j=1,\ldots,s.$
Put $U=\prod _{j=1}^s U_j^{m_j}.$
The following computations make sense for suitable small enough real
numbers $\epsilon ' >\epsilon >0.$
Write $\eta_i(v)(w_j)=f_{v,i,j}(w_j)dw_j$ on $W_j:=w_j(U_j)$ for
$i=1,\ldots,n,$
$j=1,\ldots s$ and $v\in V(\epsilon)$
The local expression $\hat{H}$ of $H$ around $q_0,$
$\hat{H}:= \Upsilon^{-1} \circ H \circ \Psi^{-1},$ is given by
$$\hat{H}:V(\epsilon)\times \xi(p_n(U))  \to V(\epsilon ')\times W(\epsilon
')$$
$$\hat{H}(v,t)=\left( v,\,
	2 \sum_{j=1}^s \sum_{h=1}^{m_j} \int_{0}^{w_{h,m_j}} f_{v,j}(w_j)dw_j
		+ C(v)- \bf{T}(v) \right)$$
where $f_{v,j}=\,^T(f_{v,1,j},\ldots,f_{v,n,j}),$
$w_{h,m_j}\equiv w_j,$
$t_{l,m_j}=\sum_{h=1}^{m_j}w_{h,m_j}^l,$
$h=1,\ldots,m_j,$ $t=(t_{1,m_j},\ldots,t_{m_j,m_j})_{j=1,\ldots,s},$
and $C(v),$ $\bf{T}(v)\in\c^n$ such that
$\pg_v(C(v))=2\varphi_v(D_0)-2\varphi_v(J_v(\infty)),$ $\pg_v({\bf
T}(v))=T(v)$ and ${\bf T}(v_0)=C(v_0).$
In order to prove that $d\hat{H}_{(v_0,0)}$ is bijective, it suffices to
see that
$$(\mbox{Jac}H_{0})|_{t=0}\ne 0$$
where $H_{0}$ is the holomorphic map defined by
$(v_0,H_{0}(t))=\hat{H}(v_0,t),$
that is to say
$$H_{0}(t)=2 \sum_{j=1}^s \sum_{h=1}^{m_j}
\int_{0}^{w_{h,m_j}}f_{v_0,j}(w_j)dw_j
$$

In the sequel we denote $f_{j}=f_{v_0,j}.$
Put $f_j(w_j)=\sum_{l=0}^{\infty} b_{j,l} w_j^l,$ $b_{j,l}\in\c^n,$
$j=1\ldots s.$
Then the Taylor series for the holomorphic map
$w_{h,m_j}\mapsto \int_{0}^{w_{h,m_j}}f_{j}(w_j)dw_j$ is
$\int_{0}^{w_{h,m_j}}f_{j}(w_j)dw_j=\sum_{l=1}^\infty a_{j,l} w_{h,m_j}^l,$
where $a_{j,l}=\frac{1}{l} b_{j,l-1},$ $l\geq 1,$ $j=1,\ldots,s.$
It is not hard to check that
$H_0(t)=2 \sum_{j=1}^s \sum_{l=1}^{m_j} a_{j,l} t_{l,m_j}+R(t),$
where the first derivatives of $R$ with respect to $t_{l,m_j}$ vanish at
$t=0,$
and so the column vectors of the Jacobian matrix of $H_{0}$ are
$\{2\,a_{l,j},l=1,\ldots,m_j,j=1\ldots,s\}.$
Reasoning by contradiction, suppose that the rows of that matrix are
linearly dependent, which is equivalent to saying that there exists a
holomorphic
1-form
$\omega_0$ in $N(v_0)$ having a zero at $z_j \in \Omega(v_0) \subset
N(v_0)$ of order at least $m_j,$ $j=1,\ldots,s.$ A direct application of
Riemann-Roch theorem gives the existence of a non-constant meromorphic
function $f$ on $N(v_0)$ having poles at $z_1,\ldots,z_s$ with order at most
$m_1,\ldots,m_s,$ respectively. In particular, $f$ has degree less than or
equal to $n.$ As $J_{v_0}(\infty)$ is not a pole of $f,$ up to adding a constant we can suppose that  $f(J_{v_0}(\infty))=0.$

On the other hand, since
$\varphi_{v_0}(D_0)-\varphi_{v_0}(J_{v_0}(\infty))=K_i(v_0)$ and
$I_{v_0}(K_i(v_0))=K_i(v_0),$ $i \in \{1,\ldots,2^{2n}\},$ then we get $\varphi_{v_0}(D_0 \cdot
\infty)-\varphi_{v_0}(J_{v_0}(D_0 \cdot \infty))=0.$ Therefore, a direct
application of Abel's theorem  gives the existence of a meromorphic
function $g$ of degree $n+1$ on $N(v_0)$ whose principal divisor coincides
with $\frac{D_0 \cdot \infty}{J_{v_0}(D_0 \cdot \infty)}.$ As $J_{v_0}$ is
an antiholomorphic involution with fixed points, it is not hard to check
that
$g \circ J_{v_0}=r/\overline{g},$ $r>0.$ Hence, up to multiplying  $g$ by the
factor $r^{-1/2},$ we can suppose that $g \circ J_{v_0}=1/\overline{g}.$
Note that $g \in {\cal F}_{v_0},$ where for any $v \in {\cal T}_n,$ ${\cal F}_v$ denotes the family of meromorphic functions $h$ of degree $n+1$ in $N(v)$ with zeroes  in  $\Omega (v)$ and satisfying
$h\circ J_{v}=1/\overline{h}.$ 
\vspace{0.3cm}
\begin{quote}{{\bf Claim:}  \em Let  $f_\lambda $ denote the meromorphic
function $\frac{1+\lambda f}{1+\overline{\lambda (f \circ J_{v_0})}},$
$\lambda \in \c.$
Then, $f_\lambda$ is not constant, for any $\lambda \in \c^*.$ Moreover,
$g_\lambda:=g f_\lambda \in {\cal F}_{v_0}$ for any $\lambda \in \c.$}
\end{quote}

Assume $f_\lambda=c,$ where  $c,$ $\lambda \in \c^*.$ Then, we infer that
$1+\lambda f=c(1+\overline{\lambda (f\circ J_{v_0})})$ and so the polar
divisor of $f,$ which is contained in $D_0,$ is invariant under $J_{v_0}.$
This is absurd because $D_0 \in Div_{n}(\Omega(v_0))$ and $\Omega(v_0) \cap
J_{v_0}(\Omega(v_0))=\emptyset.$

For the second part of the claim, first note that the principal divisor of $g_{\lambda}$ is
$[g_{\lambda}]= \frac{D_{\lambda}\cdot\infty}
{J_{v_0}(D_{\lambda})\cdot J_{v_0}(\infty)},$ where $D_{\lambda}$ is an integral divisor of degree $\leq n$ and so 
 the degree of $g_\lambda$
is $\leq n+1,$ $\lambda \in \c.$  Moreover,  $g_{\lambda}$
is not constant for any $\lambda$ (otherwise, $J_{v_0}(\infty)$ would be a zero of $1+\lambda f,$ contradicting  $f(J_{v_0}(\infty))=0$).
Let $A$ be the set $\{\lambda \in \c \;:\; g_\lambda \in {\cal F}_{v_0}\},$ and
observe that $0\in A.$ It suffices to see that $A$ is open and closed.

The openness of $A$ is  an elementary consequence of Hurwitz theorem (we
are using the fact  that the degree of $g_{\lambda}$ is at most $n+1$).
Finally, let us prove  that $A$ is closed. Let $\lambda_0 \in
\overline{A},$ and take $\{\lambda_n\}_{n \in \n} \to \lambda_0,$ where
$\{\lambda_n \;:\; n \in \n\} \subset A.$
The sequence $\{g_n:=g_{\lambda_n} \}_{n
\in \n}$ converges to $g_0:=g_{\lambda_0}$ uniformly on 
$N(v_0).$
We know  that  $g_n \circ J_{v_0}=1/\overline{g_n}$ and so the zeros of $g_n$ 
lie in
$\Omega(v_0),$ therefore, $g_n$ is holomorphic on $\Omega(v_0),$ $n\in\n$
and so the same holds for $g_0.$
Moreover, since $|g_0|=1$ on $\partial\Omega(v_0)$ and it is non constant,
the maximum principle implies that  $|g_0|<1$  on $\Omega(v_0)$ and we
infer that $g_0$ has no critical points on $\partial \Omega_{v_0}.$  As
$\partial \Omega_{v_0}$ consists of $n+1$ disjoint circles, this means that
$g_0$ takes on any complex number $\theta \in \s^1$ at least $n+1$ times.
Hence the degree of $g_0$ must be $n+1$ and $g_0 \in {\cal F}_{v_0}.$ This
concludes the proof of the claim.

\vspace{0.3cm}

To get the desired contradiction  take $P \in \partial \Omega(v_0)$ such
that $f(P) \neq 0,\;\infty,$  and choose $\lambda'= \frac{-1}{f(P)}.$
Since $J_{v_0}(P)=P,$ the meromorphic function $g_{\lambda'}$ has degree
less than $n+1,$ and so, $\lambda' \notin A=\c,$ which is absurd.

Summarizing, we have proved that
${H}|_{{\cal S}_n}: {\cal S}_n \to {\bf 0},$
${H}(v,D)=(v,0),$
is a local diffeomorphism,
where ${\bf 0}=\{(v,0)\;:\;v\in {\cal T}_n \}\subset\Jg_n$ is the null section in the Jacobian bundle.  Consequently, the projection $\vg:{\cal S}_n \to
{\cal T}_n,$ $\vg(v,D)=v,$ is a local diffeomorphism too. To finish, it suffices to check that  $\vg$ is also proper.
Indeed, take a sequence $\{(v_k,D_k)\}_{k\in\n}\subset {\cal S}_n$ such that
$\{v_k\}_{k\in\n}$ converges to a point $v_{\infty}\in {\cal T}_n.$ We can
assume that $(v_k,D_k)\in {\cal S}_n(i)$  for any $k\in\n,$ and then,
$\varphi_{v_k}(D_k)-\varphi_{v_k}(J_{v_k}(\infty))=K_i(v_k),$ $i \in \{1,\ldots,2^{2n}\}.$ Since
$I_{v_k}(K_i(v_k))=K_i(v_k),$ we get
$\varphi_{v_k}(D_k \cdot \infty)-\varphi_{v_k}(J_{v_k}(D_k \cdot
\infty)) =0.$ By Abel's theorem  there is a meromorphic function
$g_{k} \in {\cal F}_{v_k}$ with
canonical divisor $\frac{D_k \cdot \infty}{J_{v_k}(D_k \cdot
\infty)}.$ 
Let us see that $\{g_k\}_{k \in \n} \to g_\infty \in {\cal F}_{v_\infty}.$
Reflecting about all the components of $\partial \Omega(v_k),$  we can meromorphically extend $g_k$  to a {\em planar} open neighborhood $W_k$ of $\overline{\Omega(v_k))},$ $k \in \n.$ By continuity and for $k_0$ large enough, the set $W=\cap_{k \geq k_0} W_k$ is a planar neighborhood of $\overline{\Omega(v_\infty)}.$ Classical  normality criteria show that, up to taking a subsequence,
$\{g_k\}_{k\in\n}$ converges uniformly on $\overline{\Omega(v_\infty)}$ to a function $g_{\infty}$ which is meromorphic beyond $\overline{\Omega(v_\infty)}.$  It is clear that $|g_{\infty}|= 1$ on $\partial\Omega(v_{\infty}),$
$|g_{\infty}|< 1$ on $\Omega(v_{\infty})$ and
$g_{\infty}(\infty)=0.$ This proves that $g_{\infty}$ is non constant and can be extended to $N(v_\infty)$ by the Schwarz reflection
$g_{\infty}\circ J_{\infty}=1/\overline{g_{\infty}}.$ 
Since $\deg(g_k)=n+1,$ then Hurwitz theorem implies that $\deg(g_{\infty})\leq n+1.$
On the other hand, $|g_{\infty}|=1$ only on $\partial \Omega(v_{\infty}),$ and so
$g_\infty$ is injective on every boundary component of
$\Omega(v_\infty).$ Therefore, the degree of $g_\infty$ must be exactly $n+1$ and $g_\infty \in {\cal F}_{v_\infty}.$
Finally, note that $\frac{D_{\infty} \cdot \infty}{J_{v_{\infty}}(D_{\infty}\cdot \infty)}$ where $D_\infty \in \div_n$ and use Hurwitz theorem to infer that  $\{D_k\}_{k\in\n} \to D_{\infty}\in\div_n.$ Since ${\cal S}_n$ is a closed subset of $\div_n,$ we get 
$D_{\infty}\in{\cal S}_n(i),$ which proves the properness of $\vg:{\cal S}_n \to
{\cal T}_n$ and so the theorem.
\end{proof}

\section{The moduli space of maximal  graphs with singularities}
\label{sec:moduli}

The global parameterization of the space of once punctured circular domains given
in Section \ref{sec:teich} has been inspired by the theory of {\em decorated} moduli spaces of
conformal structures. For this reason, the results in the previous section connect in a natural way with the space of {\em marked} CMF graphs in $\l^3.$
We know that any CMF embedded surface is, up to an ambient isometry, a graph over
the plane $\{x_3=0\}$ with vertical limit normal vector at the end.
In the sequel we denote by $\Gg_n$ the space of CMF graphs over the plane
$\{x_3=0\}$ with vertical limit normal vector at the end and having $n+1$
singularities, $n\ge 1$.
As a consequence of Corollary \ref{co:exist}, the space $\Gg_n$ is not empty.
Let $G\in \Gg_n$  and label $F$ as its set of singularities. By definition, a {\em mark}
in $G$ is 
an ordering $\mg=(q_0,q_1,\ldots,q_n)$ of the points in $F,$ and we say that $(G,\mg)$ is a marked graph. We denote by $\Mg_n$
the space of marked graphs and call $\sg_1:\Mg_n\to\Gg_n$ and
$\sg_2:\Mg_n\to\r^{3n+4}$ the maps given by
$\sg_1(G,\mg)=G$ and $\sg_2(G,\mg)=(\mg,c),$ where $c\in\r$ is the
logarithmic growth of $G$ at the end.
By Theorem \ref{th:uniqueness}, the map
$\sg_2$ is injective.
One can use this map to endow $\Mg_n$ with an analytic structure, provided
$\sg_2(\Mg_n)$ is an open subset of $\r^{3n+4}.$ To prove this fact will be the main goal of this section.

Let $Y=(G,\mg) \in \Mg_n,$ where $\mg=(q_0,q_1,\ldots,q_n) \in \r^{3n+3},$
and orient $G$ with downward vertical limit normal vector at the end.
By Koebe's uniformization theorem, $G-F$ with the prescribed orientation is
biholomorphic to an
once punctured circular domain, where the puncture corresponds to the end
of the surface.
Hence it is not hard to see that there are {\em unique} $v \in {{\cal
T}_n}$ and conformal maximal immersion $X:\overline{\Omega(v)}-\{\infty\} \to \l^3$ such that: $G-F$ is biholomorphic to $\Omega(v)-\{\infty\}$
(in the sequel, they will be identified),
$G=X(\overline{\Omega(v)}-\{\infty\}),$ $q_0=X(\partial B_1(0)),$ and
$q_{i}=X(\partial B_{r_i(v)}(c_i(v))),$ $i=1,\ldots,n.$ 
Call $(g,\phi_3)$ the Weierstrass data for $X,$ and observe that they can be extended by Schwarz reflection to $N(v).$ Note also that $|g|<1$  on $\Omega(v).$
The behaviour of $g$ around the singularities (see Lemma \ref{lem:sing})
implies that $g:\overline{\Omega(v)} \to \{|z| \leq 1\}$  is a
$(n+1)$-branched covering, and so $g$ has exactly $n+1$ zeros
$\infty,w_1,\ldots,w_n \in \Omega(v)$ counted with multiplicity. Putting $D=w_1\cdot \ldots \cdot w_n,$ it is not hard to see that $(v,D) \in {\cal S}_n.$Indeed, the principal
divisor of $g$ in $N(v)$ is equal to $\frac{D \cdot \infty}{J_v(D \cdot
\infty)},$  and Abel's Theorem gives $\varphi_v(D \cdot
\infty)-\varphi_v(J_v(D \cdot \infty))=0.$
On the other hand, the canonical divisor of $\phi_3$ in $N(v)$ is given by
$\frac{D \cdot J_v(D)}{\infty \cdot J_v(\infty)},$ and consequently
$\varphi_v(D \cdot J_v(D))-\varphi_v(\infty \cdot J_v(\infty))=T(v).$
Therefore,
$2 (\varphi_v(D)-\varphi_v(J_v(\infty)))=T(v),$ and so, $(v,D) \in {\cal
S}_n.$

Write $\phi_3 (z) =h_3(z) \frac{dz}{z}$ on the planar domain $U(v):=
\big( \Omega(v)-\{\infty\} \big) \cup \{|z|=1\} \cup \big(
\Omega(v)^*-\{J_v(\infty)\} \big)
\subset N(v),$
where $z=\mbox{Id}|_{U(v)},$  and observe that $h_3(z) \in
\r^*$ on $|z|=1.$
\begin{definition} \label{def:cale}
With the previous notation, we call ${\cal E}: \Mg_n \to {\cal S}_n\times
\r^3\times \s^1 \times \r^*$the map given by  ${\cal E}(G,\mg)=
\big( (v,D),q_0,g(1),h_3(1) \big).$
\end{definition}
Note that the fixed orientation in the graphs is fundamental for the
definition of ${\cal E}.$ The first coordinate of  ${\cal E}$ encloses the information about the conformal structure and Weierstrass data of the marked graph, while the three last ones are simply translational, rotational and homothetical, respectively.

\begin{proposition}\label{pro:cale}
The map ${\cal E}: \Mg_n \to {\cal S}_n\times \r^3\times \s^1 \times \r^*$
is bijective.
\end{proposition}
\begin{proof}
Take $x \in {\cal S}_n,$ $x=(v,D),$ and suppose $x \in {\cal S}_n(i),$ that
is to say, $\varphi_v(D)-\varphi_v(J_v(\infty))=K_i(v).$ Since
$I_v(K_i(v))=K_i(v),$ we have
$\varphi_v(D)-\varphi_v(J_v(\infty))-\varphi_v(J_v(D))+\varphi_v(\infty)=0.$ By
Abel's theorem, there exists a meromorphic function $f$ of degree $n+1$ on $N(v)$ whose
principal divisor is equal to
$\frac{D \cdot \infty}{J_v(D \cdot \infty)}.$ Define
$$g_x:=\frac{1}{f(1)} f,$$
and observe that $g_x \circ J_v=1/\overline{g_x}$ and $g_x(1)=1.$
Moreover, the last equation and the fact $(g_x)=\frac{D \cdot \infty}{J_v(D
\cdot \infty)}$ characterize $g_x$ as  meromorphic function on $N(v).$
Since $\varphi_v(D \cdot J_v(D))-\varphi_v(J_v(\infty) \cdot \infty)=T(v),$
then there exists a meromorphic 1-form $\nu$ with canonical divisor 
$\frac{D \cdot J_v(D)}{\infty \cdot J_v(\infty)}.$
Observe that $J_v^*(\nu)=\lambda\, \overline{\nu},$ where $|\lambda|=1$
(recall that $J_v$ is an involution).
Then,  the 1-form $\phi=\frac{i}{\sqrt{\lambda}}\, \nu$ satisfies
$J_v^*(\phi)=-\overline{\phi}.$ If we put $\phi(z)=h(z)\, \frac{dz}{z},$
$z\in U(v),$ we infer that $h(z) \in \r^*,$ $|z|=1.$
Then, define $$\phi_3(x):=\frac{1}{h(1)}\, \phi,$$ and observe that the
equations
$(\phi_3(x))=\frac{D \cdot J_v(D)}{\infty \cdot J_v(\infty)}$ and
$h_3(1)=1$ characterize $\phi_3(x)$ as meromorphic 1-form on $N(v).$

With this notation, given $\big( x,q_0,\theta,r \big) \in {\cal S}_n\times
\r^3\times \s^1 \times \r^*,$ define
\begin{equation} \label{eq:datos}
g_x(\theta)= \theta g_x, \quad \phi_3(x,r)= r\, \phi_3(x)
\end{equation}
and define $\phi_1(x,\theta,r)$ and $\phi_2(x,\theta,r)$ in the obvious way.
As $g_x(\theta)$ is holomorphic on $\Omega(v)$ and $|g|=1$ on
$\partial \Omega(v),$ the maximum principle implies that $|g|<1$ on
$\Omega(v),$ and so by Theorem \ref{th:exist} the map $X_x(q_0,\theta,r):\overline{\Omega(v)}-\{\infty\} \to \l^3,$
$$X_x(q_0,\theta,r)(z):= q_0+\mbox{Real} \int_1^z
\big(\phi_1(x,\theta,r),\phi_2(x,\theta,r),\phi_3(x,r) \big)$$
provides a CMF graph $G_x(q_0,\theta,r):=X_x(q_0,\theta,r)\Big(\overline{\Omega(v)}-\{\infty\}\Big)
\in\Gg_n.$
 
Defining the mark $\mg_x(q_0,\theta,r)$ by
$q_0=X_x(q_0,\theta,r)(\partial
B_{0}(1)),$   $q_{i}=X_x(q_0,\theta,r)(\partial
B_{r_i(v)}(c_i(v))),$ $i=1,\ldots,n,$ it is now clear that
${\cal
E}^{-1}(x,q_0,\theta,r)=\{(G_x(q_0,\theta,r),\mg_x(q_0,\theta,r))\},$ and
so, ${\cal E}$ is bijective.

\end{proof}

Label ${\cal P}_{n+1}$ as the symmetric group of permutations of order $n+1.$ We
denote by $\mu:{\cal P}_{n+1} \times \Mg_n \to \Mg_n,$  the natural action
$\mu(\tau,(G,\mg)):=(G,\tau(\mg)).$

The following theorem will show that $\Mg_n$ and $\Gg_n$ are analytic
manifolds of dimension $3n+4.$ We first need the following lemma:

\begin{lemma} \label{lem:porfin}
Given $v \in {{\cal T}_n},$ there exists a holomorphic 1-form $\omega_0$ in
$N(v)$ having $2n-2$ distinct zeroes, none of them contained in $\partial \Omega(v),$ and satisfying $J_v^*
(\omega_0)=\overline{\omega_0}.$
\end{lemma}
\begin{proof}
Let $\omega$ be a holomorphic non-zero 1-form on $N(v)$ satisfying
$J_v^*(\omega)=\overline{\omega}.$ In addition, we can choose $\omega$ such
that it does not vanish on $\partial\Omega(v).$
Indeed, let $h$ be the unique harmonic function on $\overline{\Omega(v)}$
satisfying $h|_{a_{j}(v)}=\delta_{j,0},$ $j=0,1,\ldots,n,$ where $\delta_{h,k}$ is the Kronecker symbol.
By the maximum principle,  $h^{-1}(0)=\cup_{j=1}^n a_j(v)$ and
$h^{-1}(1)=a_{0}(v).$ Since these level sets are regular, $h$ has no
critical point on
$\partial \Omega(v).$ The 1-form  $\omega= \partial_z h$ can be extended to
$N(v)$ by Schwarz reflection, and satisfies the desired properties.

Let $D=P_1^{n_1} \cdot \ldots  \cdot P_s^{n_s} \cdot J_v(P_1)^{n_1} \cdot
\ldots  \cdot J_v(P_s)^{n_s}$ be the canonical divisor
associated to $\omega.$

Since the degree of $D$ is $2n-2,$ it is not hard to deduce from
Riemann-Roch theorem that the complex linear space ${\cal M}(D)$ of
meromorphic functions on $N(v)$ having poles only at the points $P_i,$
$J_v(P_i)$ with order at most $n_i,$ $i=1,\ldots,s,$ has dimension $n.$

Note that there is a function $f \in {\cal M}(D)$ of degree $2n-2.$ Indeed,
otherwise the maximum degree $d_0$ of functions in ${\cal M}(D)$ would be
less than $2n-2,$ and so, there would exist an entire divisor $D'$ of
degree $d_0$ such that $D \geq D'$ and ${\cal M}(D')={\cal M}(D),$ where
${\cal M}(D')$ is
defined in a similar way.    Given $Q \in N(v)$ such that $ D/D' \geq Q,$ then,
with obvious notations, ${\cal M}(D') \subset {\cal M}(D/Q)  \subset {\cal
M}(D).$
We infer that ${\cal M}(D') = {\cal M}(D/Q)  = {\cal M}(D)$ and  so ${\cal
M}(D/Q)$ has dimension $n.$ It then follows from Riemann-Roch theorem that
the
complex linear space of holomorphic 1-forms $\omega'$  satisfying
$(\omega') \geq D/Q$ has dimension $2.$ Take a holomorphic 1-form
$\omega' \in H$ linearly independent from $\omega.$ The quotient
$\omega/\omega'$ then defines a meromorphic function on $N(v)$ of degree
1, which is absurd.

Hence, we can take $f\in{\cal M}(D)$ with $\deg f=2n-2.$ Since the polar divisor of
$f$ is $D$
and $J_v(D)=D,$  we can find $\lambda$ and $\mu \in \r$ such that
$f_0:=\lambda (f+\overline{f \circ J_v})+i \mu (f-\overline{f \circ J_v})$
has the same polar divisor as $f,$ moreover, one has $f_0 \circ
J_v=\overline{f_0}.$
Since poles of $f_0$ do not lie in $\partial \Omega(v),$ $f_0
|_{\partial\Omega(v)}$ is bounded. Therefore, there exists a large enough
real number $r$ such that $f_0+r$ has only simple zeroes, and these zeroes
are not contained in $\partial \Omega(v).$ The 1-form $\omega_0:=(f_0+r)
\omega$ satisfies the properties in the statement of the lemma.
\end{proof}
In Subsection \ref{subsec:spin},  the smoothness of $\hat{T}$ was derived from the one of the Riemann's constants vector.  However, this fact can be also deduced from the previous lemma, at least for the first order derivatives. Indeed, let $v_0 \in {{\cal T}_n}$ and consider a
holomorphic 1-form $w(v_0)$ on $N(v_0)$ with simple zeroes (for instance, the one constructed in the lemma). Put $w(v_0)=
\sum_{j=1}^n \lambda_j(0) \eta_j(v_0)$ and note that $\lambda_j(0)\in\r.$ Label $\lambda(0)= (\lambda_1(0),\ldots,\lambda_n(0)),$ and without loss of generality, suppose $\lambda_1(0)=1.$ Define
$w_\lambda(v)= \sum_{j=1}^n \lambda_j \eta_j(v),$ $v \in {\cal T}_n,$ $\lambda=(1,\lambda_2,\ldots,\lambda_n)\in \{1\} \times \r^{n-1}$ and observe
that $J_v^*(\omega_\lambda(v))=\overline{\omega_\lambda(v)}.$
A direct application of Hurwitz theorem implies that $w_\lambda(v)$ has also simple
zeros
which are not contained in $\partial\Omega(v),$ for $v$  and $\lambda$ close enough to $v_0$ and $\lambda(0),$ respectively. 
Moreover, the Implicit Function theorem and Corollary \ref{co:smooth} give that the zeros of $w_\lambda(v)$
depend at least ${\cal C}^1$  on $(v,\lambda)$ (in fact smoothly, see the Remark below). Taking into account the smoothness of  
$\varphi:\div_{n-1} \to {\cal J}_n$ and 
${\cal I}:{\cal J}_n\to{\cal J}_n,$ we infer that $\hat{T}:{{\cal T}_n} \to {\cal J}_n$ is at least ${\cal C}^1.$ 

\begin{remark} \label{re:porfin}
Since $\hat{T}$ and $\varphi$ are  smooth, then the zeroes of $w_\lambda(v),$ for $(\lambda,v)$ in a neighbourhood of $(\lambda(0),v_0),$ locally  define the following  $4n-2$-dimensional {\em smooth} submanifold of $\div_{n-1}:$ $$\{(v,D) \in \div_{n-1} \;:\; \varphi(v,D)+\varphi(v,J_v(D))=\hat{T}(v)\}.$$ 
\end{remark}

\begin{proposition} \label{pro:ruso}
Endowing $\Mg_n$ with the unique
differentiable structure making ${\cal E}$ a diffeomorphism, the maps
${\cal S}_n \to {\cal C}_n,$ $ x \mapsto g_x,$ and ${\cal S}_n \to {\cal H}_n,$ $x \mapsto \phi_3(x),$  are smooth with $2$-regularity and $1$-regularity, respectively. 
\end{proposition}
\begin{proof}
Indeed, take $x_0=(v_0,D_0) \in {\cal S}_n.$
From Theorem \ref{th:submanifold},
there exists an open ball $V(\epsilon)$ in ${\cal T}_n$  centered at $v_0$
of radius $\epsilon >0$ and a local diffeomorphism ${V}(\epsilon) \to
{\cal S}_n,$ $v \mapsto (v,D(v)),$ where $D(v_0)=D_0.$ For simplicity, we
write $x(v):=(v,D(v)),$ $v \in { V}(\epsilon).$

Therefore,  the map ${V}(\epsilon) \to \div_{n+1},$ $v \to (v,\infty \cdot
D(v))$ is smooth, and since $\tau:\div_{n+1} \to {\cal H}_n$ is also smooth with $1$-regularity
(see Corollary \ref{co:smooth}), the same holds for the map ${V}(\epsilon)
\to {\cal H}_n,$ $v \mapsto \tau_{v}:=\tau_{\infty \cdot D(v)}(v).$

Take a smooth deformation of $N(v_0),$ $\{F_v:N(v_0)\to N(v)\;:\;v \in
V(\epsilon)\}.$  Let
$B(v_0)=\{a_1(v_0),\ldots,a_n(v_0),b_1(v_0),\ldots,b_n(v_0)\}$ be the
canonical homology basis on $\Omega(v_0)$ defined as in Section
\ref{sec:teich}. In what follows, we deal with any representative curves $a_j(v_0),$ $b_j(v_0),$ $j=1,\ldots,n,$ 
of these homology classes for which $N(v_0)-\cup_{j=1}^n (a_j(v_0) \cup b_j(v_0))$ is simply connected and contains  the points in $\infty \cdot D_0.$
For small enough $\epsilon,$ the curves
$a_j(v):=F_v(a_j(v_0)),$ $b_j(v_0):=
F_v(b_j(v_0))$ do not pass also through the points in $ \infty \cdot
D(v),$ $v \in {V}(\epsilon),$ $j=1,\ldots,n.$

y Abel's theorem, and for $z \in N(v)-\cup_{j=1}^n (a_j(v) \cup b_j(v)):$
$$g_{x(v)}(z)=\mbox{Exp} \big( \int_1^z (\tau_{v}+ \sum_{j=1}^n m_j(v)
\eta_j(v))  \big).$$ In this expresion, the integration paths lie in $\big(N(v)-\cup_{j=1}^n (a_j(v) \cup
b_j(v))\big) \cup\{1\},$
and $m_j(v) \in \z$ are integer numbers determined by the equation:
$$\widetilde{\varphi_v}(\infty \cdot
D(v))-\widetilde{\varphi_v}(J_v(\infty) \cdot J_v(D(v)))=\sum_{j=1}^n
m_j(v) \pi^j(v),$$ where $\widetilde{\varphi_v}$ is the branch of
$\varphi_v$ on $N(v)-\cup_{j=1}^n (a_j(v) \cup b_j(v))$ vanishing at $1.$

Since $m_j(v)$ depend continuously on $v,$ then $m_j(v)=m_j \in \z$ and
so, by Corollary \ref{co:smooth}, $g_{x(v)}$ depends smoothly on $v$ with $2$-regularity.

We have to obtain the analogous result for the map
${\cal V}(\epsilon) \to {\cal H}_n,$ $v \to \phi_3(x(v)).$
Take the holomorphic 1-form $\omega_0$ on $N(v_0)$ given in Lemma
\ref{lem:porfin}, write $\nu(v_0):=\omega_0=\sum_{j=1}^n \lambda_j
\eta_j(v_0),$ where $\lambda_j \in \r,$ and define $\nu(v):=\sum_{j=1}^n
\lambda_j \eta_j(v).$
Since the map $v\mapsto \nu(v)$ is smooth with $1$-regularity (see Corollary \ref{co:smooth})
it suffices to prove that $v\mapsto \frac{\phi_3(x(v))}{\nu(v)}$ is smooth with $2$-regularity.
By Hurwitz's Theorem and the implicit function theorem, $\nu(v)$
satisfies also the thesis in Lemma \ref{lem:porfin}, for small enough
$\epsilon.$ Moreover, as explained during the proof of Lemma \ref{lem:porfin}, the map $V(\epsilon) \to \div_{2n-2},$ $v \mapsto (v,(\nu(v)))$ is at least ${\cal C}^1$ (in fact smooth by Remark \ref{re:porfin}), where as usually $(\nu(v))$ is the canonical divisor associated to $\nu(v).$
Hence,  writing $(\nu(v))=E(v) \cdot J_v(E(v)),$ the
map $V(\epsilon) \to \div_{n-1},$ $v \mapsto E(v)$ is also smooth, and 
therefore, the same holds for  ${ V}(\epsilon) \to \div_{n,n},$ $v \mapsto
(v,\infty\cdot E(v),D(v)).$ We infer from Corollary \ref{co:smooth} that the  map ${V}(\epsilon) \to {\cal H}_n,$ $v \mapsto \kappa_{v}:=\kappa_{
\infty \cdot E(v),D(v)}(v)$ is smooth with 1-regularity.
Reasoning as above, the map
$$f_{x(v)}(z)= \mbox{Exp} \big( \int_1^z (\kappa_{v}+ \sum_{j=1}^n n_j
\eta_j(v))  \big),$$
is a well defined meromorphic function on $N(v),$
for suitable integer numbers $n_j$  not depending on $v$ and
${\cal V}(\epsilon) \to {\cal C}_n,$ $v \mapsto f_{x(v)},$
is smooth with $2$-regularity. The principal divisor associated to $f_{x(v)}$ is given by
$(f_{x(v)})=\frac{D(v) \cdot J_v(D(v))}{\infty \cdot E(v) \cdot J_v(\infty)
\cdot J_v(E(v))}.$
Therefore, if we write $\nu(v)=h_v(z) \frac{dz}{z}$ on $U(v),$
we infer that $\frac{\phi_3(x(v))}{\nu(v)}=\frac{1}{h_v(1)} f_{x(v)},$ and
so $v\mapsto \phi_3(x(v))$ is smooth with $1$-regularity. This concludes the proof.
\end{proof}
\begin{remark}
Equation (\ref{eq:datos}) gives that the maximal immersion
$X_x(q_0,\theta,r)$ determining the marked graph
$(G_x(q_0,\theta,r),\mg_x(q_0,\theta,r))$ depends smoothly on $(x,q_0,\theta,r)$ with $2$-regularity (the notion of differentiability with $k$-regularity for immersions is defined in the obvious way).
\end{remark}

\begin{theorem}[Main theorem] \label{th:structure}
The set $\sg_2(\Mg_n) \subset \r^{3n+4}$ is open and hence the one to one map
$\sg_2:\Mg_n \to \sg_2(\Mg_n)$ provides a global system of analytic coordinates on
$\Mg_n.$
Moreover, the action $\mu$ is discontinuous and hence the orbit space,
naturally identified to $\Gg_n,$ has a unique analytic structure making
$\sg_1$ an analytic covering of  $(n+1)!$ sheets.
\end{theorem}

\begin{proof}
Let $x\in{\cal S}_n$ and denote by $(G_x,\mg_x)={\cal
E}^{-1}(x,0,1,1)\in\Gg_n.$
Call $X_x$  the associated maximal immersion and label as $(g_x,\phi_3(x))$
its Weierstrass data.

Observe that  ${\cal E}^{-1} (x \times \r^3 \times \s^1 \times ]0,+\infty[)$
consists of all the marked graphs which differ from $(G_x,\mg_x)$
by  translations, rotations about a vertical axis (i.e., parallel to the $x_3$-axis) and homotheties.
Since $X_x(q_0,\theta,r)$ depends smoothly on $(x,q_0,\theta,r)$ with $2$-regularity, the map $\sg_2:\Mg_n \to \r^{3n+4}$ is smooth
By the injectivity of $\sg_2$ (see  Theorem \ref{th:uniqueness}) and the domain invariance theorem, ${\sg}_2(\Mg_n)$ is open in $\r^{3n +4}$
and hence it is an analytic manifold of dimension $3n+4.$ We can then
endow $\Mg_n$ with the unique analytic structure making ${\sg}_2: {\Mg}_n
\to {\sg}_2({\Mg}_n)$  an analytic diffeomorphism. 

To conclude,  it remains to check that the action $\mu$ is discontinuous. Indeed, let $\tau:\Mg_n
\to \Mg_n$ denote the diffeomorphism given by $\tau(G,\mg)=(G,\tau(\mg)),$
$\tau \in {\cal P}_{n+1}.$ Let $(G_0,\mg_0) \in \Mg_n$ and write
$\mg_0=(q_0,q_1,\ldots,q_n) \in \r^{3n+3}.$ Take a neighborhood $U_j$ of
$q_j$ in $\r^3,$ $j=0,1,\ldots,n,$ such that $U_i \cap U_j=\emptyset,$ $i
\neq j,$ and call ${\cal U}=\prod_{j=0}^n U_j.$  Then, it is clear that
$\tau(\sg_2^{-1}({\cal U}\times \r)) \cap \sg_2^{-1}({\cal U}\times \r)=\emptyset,$ for any $\tau \in
{\cal P}_{n+1}-\{\mbox{Id}\},$ which proves the discontinuity of $\mu$ and
concludes the proof.
\end{proof}

During the proof of Theorem \ref{th:structure} we have shown that the
topological structures induced by
${\cal E}$ and $\sg_2$ on $\Mg_n$ are the same. The following theorem proves
that this topology coincides with the one of the uniform convergence of graphs on compact subsets of $\{x_3=0\}.$

\begin{theorem}
Let $\{G_k\}_{k \in \n}$ be a sequence in $\Gg_n,$ and $G_0 \in \Gg_n.$ 

Then $\{G_k)\}_{k \in \n} \to G_0$ in the topology of $\Gg_n$ if and only if 
$\{G_k\}_{k \in \n}$ converges to $G_0$ uniformly on compact subsets of $\{x_3=0\}.$   
\end{theorem}
\begin{proof} 
Suppose $\{G_k\}_{k \in \n} \to G_0 \in \Gg_n$ in the topology of $\Gg_n,$ and 
choose marks in such a way that $\{(G_k,\mg_k)\}_{k \in \n}$  converges  to  $(G_0,\mg_0) $ in $\Mg_n.$

Write ${\cal E}((\Gg_k,\mg_k))=(x_k,q_0(k),\theta_k,r_k),$ $X_k=X_{x_k}(q_0(k),\theta_k,r_k)$ and $x_k=(v_k,D_k) \in {\cal S}_n,$ 
$k \in \n \cup \{0\}.$ 

Let $W$ be any compact domain in $\r^2 \equiv \{x_3=0\}$ containing the singularities in $\mg_0$ as interior points, and let $W_k$ denote the compact set $X_k^{-1}(W \times \r) \subset \overline{\Omega}(v_k)-\{\infty\},$ $k \in \n \cup \{0\}.$ 

Since $\{(x_k,q_0(k),\theta_k,r_k)\} \to (x_0,q_0(0),\theta_0,r_0)$ and $X_{x}(q_0,\theta,r)$ depends smoothly on $(x,q_0,\theta,r)$ with 2-regularity,  it is not hard to check that $\lim_{z \to \infty} ||X_k(z)||=+\infty$ uniformly in $k.$ 
In addition, the domains $W_k$ are uniformly bounded in $\c,$ $\{W_k\}_{k \in \n} \to W_0$ in the Hausdorff distance and $X_k$ converges uniformly on $W_0$ to $X_0.$ In the last statement we have used that $X_k$ can be reflect analytically about the circles in $\partial {\Omega}(v_k),$ and so all the immersions $X_k,$ $k$ large enough, are well defined in a universal neighborhood of $W_0$ in $\c.$ It is then obvious that the function $u_k: \r^2 \to \r$ defining the graph $G_k$ converges uniformly over $W$ to the function $u_0: \r^2 \to \r$ defining $G_0$ (furthermore, $\{v_k\}_{k \in \n} \to v_0$ implies that
$\{\mg_k\}_{k \in \n} \to \mg_0$). Since $W$ can be as larger as we want, $\{u_k\}_{k \in \n} \to u_0$ uniformly on compact subsets of $\r^2.$   

Assume now that $\{G_k\}_{k \in \n}$ converges to $G_0$ uniformly on compact subsets of $\{x_3=0\},$ and as above, denote by $u_k:\r^2 \to \r$  the function defining the graph $G_k,$  $k \in \n \cup \{0\}.$

Let us show that  singular points of $G_0$ are  limits of sequences of singular points of graphs $G_k,$ $k \in \n.$ 
Indeed, let $p_0=(y_0,u_0(y_0)) \in G_0$ be a singular point, and without loss of generality, suppose that $p_0$ is a downward pointing conelike singularity. By Lemma \ref{lem:sing}, there exists $\epsilon>0$ small enough such that $u_0^{-1}(\{x_3 \leq u_0(y_0)+\epsilon\})$ contains a compact component $C_0(\epsilon)$ with regular boundary and containing $y_0$ as the unique (interior) singular point. 
Since $\{u_k\}_{k \in \n} \to u_0$ uniformly on compact subsets, $u_k^{-1}(\{x_3 \leq u_0(y_0)+\epsilon\})$ 
must contain a compact component $C_k(\epsilon)$ containing $y_0$ as well, $k$ large enough. Furthermore, $\{C_k(\epsilon)\} \to C_0(\epsilon)$ in the Hausdorff sense, and by the maximum principle $C_k(\epsilon)$ must contain at least an interior singular point $y_k$ of $u_k,$ $k$ large enough. Since $C_0(\epsilon)$ converges to $\{y_0\}$ as $\epsilon \to 0,$ we deduce that $\{p_k:=(y_k,u_k(y_k))\}_{k \to \infty} \to p_0.$

As a consequence, there exist marked graphs $(G_k,\mg_k) \in \Mg_n,$ $k \in \n \cup \{0\},$ such that $\{\mg_k\}_{k \in \n} \to \mg_0.$

Call $c_k$ the logarithmic growth of $G_k,$ $k \in \n \cup \{0\},$ and let us see that $\{c_k\} \to c_0.$ 
Indeed, let $\gamma$ be a circle in $\r^2$ containing all the singular points of $u_0$ in its interior, and let $A$ denote a closed tubular neighbourhood 
of $\gamma$ in $\r^2$ not containing any singular point of $u_0.$ It is well known that the function $u_k-u_0$ is solution of a uniformly elliptic linear  equation $L_k(u_k-u_0)=0$ over $A,$ $k$ large enough. Moreover, the fact that the functions $\frac{1}{1-|\nabla u_k|},$ $k \in \n,$ are uniformly bounded on $A$ (see \cite{Bar-Sim}) guarantee that the coefficients of operators $L_k,$ $k \in \n,$ are uniformly bounded too. Therefore, since $\{u_k\}_{k \in \n} \to u_0$ uniformly on $A,$   the classical Schauder estimates (\cite{gilbarg} p. 93) imply that $\{u_k\}_{k \in \n} \to u_0$ in the $C^2$ norm on $A.$ In particular, $$\{ \int_\gamma \nu_k(s_k) ds_k  \}_{k \in \n} \to \int_\gamma \nu_0(s_0) ds_0,$$ 
where $\nu_k$ and $s_k$ are the conormal vector and the arc-length parameter along $\gamma$  in $G_k,$ respectively, for any $k \in \n \cup \{0\}.$ 
Since $\int_\gamma \nu_k(s_k) ds_k=2 \pi (0,0,c_k)$ and $\int_\gamma \nu_0(s_0) ds_0=2 \pi (0,0,c_0),$ we infer that $\{c_k\} \to c_0.$ 

Since $\sg_2:\Mg_n \to \sg_2(\Mg_n) \subset \r^{3n+4}$ is an homeomorphism,  $\{(G_k,\mg_k)\}_{k \in \n} \to (G_0,\mg_0)$ in the manifold $\Mg_n,$ and so, 
$\{G_k\}_{k \in \n} \to G_0$ in the manifold $\Gg_n.$ This concludes the proof.
\end{proof}

It is also interesting to 
notice that $\Mg_n$ is far from being connected, because the sets ${\cal
S}_n(i),$
$i=1,\ldots,2^{2 n},$ are pairwise disjoint and open in ${\cal S}_n.$
Although the discussion could be a little more involved, it is possible to
see that the
same holds for $\Gg_n.$
On the other hand, it is natural to ask whether  the differentiable
structure in $\Mg_n$ induced by $\sg_2$ coincides with the one  induced by
${\cal E}$ or not. However, this fact is not relevant for the subsequent
analysis, and for this reason we have prefered to leave it as an open question.
In the sequel, and up to an explicit mention of the contrary, the underlying structure in $\Mg_n$ will be the analytic one induced by $\sg_2.$
\vspace{0.2cm}

There is a natural connection between the n-spinorial bundle ${\cal S}_n$
and a suitable quotient of $\Mg_n.$ Next lemma and corollary will be
devoted to give a detailed explanation of this fact. We first fix the following notations.

Let ${\cal R}$ denote the group of Lorentzian similarities generated by
translations, rotations about a vertical axis, symmetries with respect to
a horizontal plane and homotheties in $\l^3.$

Denote by $\mbox{proj}_0:{\cal S}_n \times \r^3 \times \s^1 \times \r^* \to
{\cal S}_n$ the natural projection $\mbox{proj}_0(x,q,\theta,r)=x.$

Let $\pi:\r^3 \to \{x_3=0\}$ denote the orthogonal projection, and call $\Delta^{3n+4}$ the open subset of $\r^{3n+3} \times \r \equiv \r^{3n+4}$ given by $$\Delta^{3n+4}=\{((q_0,q_1,\ldots,q_n),c)\in \r^{3n+4} \;:\; \pi(q_i) \neq \pi(q_j),\; i \neq j\}.$$
Let $\Cg_n$be the quotient of $\Delta^{3n+4}$
under the congruence relation: $((q_0,q_1, \ldots,q_n),c)
\sim ((q'_0,q'_1, \ldots,q'_n),c')$ if there exists a Lorentzian similarity
$\hat{R} \in {\cal R}$  such that $\hat{R}(q_0,q_1, \ldots,q_n)=(q'_0,q'_1,
\ldots,q'_n)$ and $\vec{R}(0,0,c)=(0,0,c'),$ where $\vec{R}$ is the linear
similarity associated to $\hat{R}.$ 
Label $\mbox{proj}_1:\Delta^{3n+4} \to
\frac{\Delta^{3n +4}}{\sim}$ as the natural projection, and denote by
$H_i=\{ ((0,q_1, \ldots,q_n),c) \;:\; x_2(q_1)=0,\; |q_1|=1,\;
x_1(q_1),x_3(q_i)>0\},$ $i=1,\ldots,n,$ $H_{n+1}=\{ ((0,q_1, \ldots,q_n),c) \;:\;
x_2(q_1)=0,\; |q_1|=1,\;x_1(q_1)>0 ,\; c>0\},$ 
$H_0=
\{((0,q_1,\ldots,q_n),0) \;:\; x_2(q_1)=0,\; |q_1|=1,\; x_1(q_1)>0,\;x_3(q_i)=0, \;
i=1,\ldots,n \},$ where $|q_1|$ here refers to the Euclidean norm of $q_1.$ 

It is clear that $\cup_{j=1}^{n+1} \mbox{proj}_1(H_j)=\frac{\Delta^{3n
+4}}{\sim}-\mbox{proj}_1(H_0).$ The map $\mbox{proj}_1|_{H_i}:H_i \to \mbox{proj}_1(H_i)$ is bijective and
provides analytic coordinates on the open set $\mbox{proj}_1(H_i) \subset
\frac{\Delta^{3n +4}}{\sim}.$ Hence, $$\Cg_n:=\frac{\Delta^{3n
+4}}{\sim}-\mbox{proj}_1(H_0)$$ is an analytic manifold.

Moreover,  it is
easy to check that $\sg_2(\Mg_n)\subset \Delta^{3n+4}-\widetilde{H}_0,$  where
$\widetilde{H}_0=\mbox{proj}_1^{-1}(\mbox{proj}_1(H_0)).$ Otherwise, there would be $(G,\mg) \in \Mg_n$ such that $\sg_2((G,\mg))\in H_0,$ and so,  $G$ would be a graph asymptotic to a horizontal plane and with singularities  in   $\{x_3=0\}.$ Remark \ref{re:uniqueness} would imply that $G=\{x_3=0\},$
which is absurd.

As a consequence, 
$\mbox{proj}_1(\sg_2 (\Mg_n)) \subset \Cg_n.$
\vspace{0.2cm}

Two marked graphs $(G_1,\mg_1),$  $(G_2,\mg_2)\in \Mg_n$ are said to be
congruent if there exists a Lorentzian similarity $\hat{R}\in {\cal R}$
such that $\hat{R}(G_1)=G_2,$ $\hat{R}(\mg_1)=\mg_2.$ Call $\hat{\Mg}_n$
the quotient space of $\Mg_n$ under the congruence relation and  denote by
$\mbox{proj}_2:\Mg_n \to \hat{\Mg}_n$ the natural projection.

Likewise $G_1,$  $G_2 \in \Gg_n$ are said to be congruent if there exists a
Lorentzian similarity $\hat{R} \in {\cal R}$ such that $\hat{R}(G_1)=G_2.$
We call $\hat{\Gg}_n$ the quotient space of $\Gg_n$ under the congruence
relation, and label $\mbox{proj}_3:\Gg_n \to \hat{\Gg}_n$ the natural
projection.

The connection between the above projections is given in the following lemma.
\begin{lemma}  \label{lem:proj}
There exists unique maps $\hat{\sg_2}:\hat{\Mg}_n \to \Cg_n,$
$\hat{\sg_1}:\hat{\Mg}_n \to \hat{\Gg}_n$
and $\hat{{\cal E}}: \hat{\Mg}_n \to {\cal S}_n$ making commutative the
following diagrams:
$$\begin{CD}
\r^{3n+4}-\widetilde{H}_0 @<\sg_2<< \Mg_n @ >{\cal E}>> {\cal S}_n\times
\r^3\times \s^1 \times  \r^* \\
@VV proj_1 V @ VV proj_2 V@ VV proj_0 V \\
{\Cg}_n @< \hat{\sg}_2<< \hat{\Mg}_n @> \hat{\cal E}>> {\cal S}_n \\
\end {CD}$$

$$\begin{CD}
\Mg_n @>\sg_1 >>\Gg_n  \\
@V proj_2 VV @ VV proj_3 V \\
\hat{\Mg}_n @>\hat{\sg}_1>>\hat{\Gg}_n \\
\end {CD}$$
Moreover, $\hat{\sg}_2$ is injective and $\hat{{\cal E}}$ is bijective.
\end{lemma}
\begin{proof}
It is natural to define $$\hat{{\cal E}}: \hat{\Mg}_n \to {\cal S}_n, \quad
\hat{{\cal E}} (\mbox{proj}_2(G,\mg))=\mbox{proj}_0({\cal E}(G,\mg)).$$ To
prove that it is well defined and bijective, it suffices to check that, given 
$(G_1,\mg_1),$ $(G_2,\mg_2) \in \Mg_n,$
$\mbox{proj}_0({\cal E}((G_1,\mg_1)))=\mbox{proj}_0({\cal E}((G_2,\mg_2)))$
if and only if  there exists $\hat{R} \in {\cal R}$ such that
$\hat{R}(G_1)=G_2,$ $\hat{R}(\mg_1) =\mg_2.$ 
Indeed, write  ${\cal E}((G_j,\mg_j))=(x_j,q_{0,j},\theta_j,r_j),$  where
$x_j=(v_j,D_j) \in {\cal S}_n,$ $j=1,2.$ For simplicity,  put
$X_j=X_{x_j}(q_{0,j},\theta_j,r_j),$ $j=1,2,$ and identify $G_j \equiv
{\overline{\Omega}(v_j)} -\{\infty\},$ $j=1,2.$
Assume there exists $\hat{R} \in {\cal R}$ such that $\hat{R}(G_1)=G_2,$
$\hat{R}(\mg_1) =\mg_2,$
and write by $R:\overline{\Omega(v_1)} \to \overline{\Omega(v_2)}$ the
natural conformal transformation induced by $\hat{R}.$
If we orient $G_1$ and $G_2$ with downward limit normal vector at the end,
any translation, rotation about a vertical axis,
homothety or symmetry with respect to a horizontal plane preserves the
orientation of the graphs, and so, the same holds for $\hat{R}.$
Therefore, $R$ is a biholomophism (that is to say, a Möbius
transformation). Since $\hat{R}$ also preserves the end and the singular points, we infer  that
$R(\infty)=\infty,$ $R(\{|z|=1\})=\{|z|=1\},$  and
$R(\partial B_{r_j(v_1)}(c_j(v_1)))=\partial B_{r_j(v_2)}(c_j(v_2)),$ for any $j.$ This implies that $v_1=v_2$ and  $R$ is the identity
map on $\overline{\Omega}(v_1).$

On the other hand, $g_{x_2} (\theta_2)  \circ R=L \circ g_{x_1}(\theta_1),$
where $L:\overline{\c} \to \overline{\c}$ is
the conformal transformation induced by $\hat{R}.$ In addition, $L$ fixes
the origin, because $R$ preserves the end and the limit normal
vector at the
ends points downward.Therefore $g_{x_2}(\theta_2)=g_{x_2}(\theta_2) \circ R= \theta  g_{x_1}(\theta_1),$ where $|\theta|=1,$ and hence $D_1=D_2.$ Now it is straightforward to check that $
x_1=x_2.$

For the converse, take $(G_1,\mg_1) \in \Mg_n,$ call ${\cal
E}(G_1,\mg_1)=(x_1,q_1,\theta_1,r_1) ,$ and observe that
the set $(\mbox{proj}_0 \circ {\cal E}
)^{-1}(x_1)=\{\big(X_{x_1}(q,\theta,r)((\overline{\Omega}(v_1)-\{\infty\}),\mg_{
x_1}(q,\theta,r))\big) \;:\;
(q,\theta,r)
\in \r^3 \times \s^1 \times \r^*\}.$ But this set consists of  marked
graphs in $\Mg_n$ differing from $(G_1,\mg_1)$ by  ambient similarities
$\hat{R} \in
{\cal R}$ preserving the mark.
Finally, observe that $\hat{\sg_2}:\hat{\Mg}_n \to \Cg_n,$
$\hat{\sg_2}([(G,\mg)]):=[(\mg,c)]$ is well defined, and use Theorem
\ref{th:uniqueness} to show that it is injective.
The map $\hat{\sg_1}:\hat{\Mg}_n \to \hat{\Gg}_n$ is given by
$\hat{\sg_1}([(G,\mg)])=[G],$ and the commutativity of the diagrams is obvious.
\end{proof}

From Lemma \ref{lem:proj}, $\hat{\Mg}_n$ can be endowed with the
differentiable structure making $\hat{{\cal E}},$ $\mbox{proj}_2$ and
$\hat{\sg}_2$ a diffeomorphism, a submersion and a smooth map, respectively
(provided that $\Mg_n$ is endowed with the differentiable structure induced
by ${\cal E}$).   However following the aim of Theorem \ref{th:structure}
we are interested in viewing $\hat{\Mg}_n$ as an analytic manifold. This
is the main goal of the following corollary.
\begin{corollary}
The set $\hat{\sg}_2(\hat{\Mg}_n) \subset \Cg_n$  is open, and  the
bijective map $\hat{\sg}_2:\hat{\Mg}_n \to \hat{\sg}_2(\hat{\Mg}_n)$
provides a unique analytic structure in $\hat{\Mg}_n$ making
$\mbox{proj}_2:\Mg_n \to \hat{\Mg}_n$ an analytic submersion.
Moreover, if $\hat{\Sigma}_1,$ $\hat{\Sigma}_2$ are two distinct connected
components of $\hat{\Mg}_n,$ then either
$\hat{\sg_1}(\hat{\Sigma}_1) \cap \hat{\sg_1}(\hat{\Sigma_2})=\emptyset$
or $\hat{\sg_1}(\hat{\Sigma}_1) = \hat{\sg_1}(\hat{\Sigma}_2).$ In the
second case, $\hat{\Sigma}_1$ and $\hat{\Sigma}_2$ are analytically
diffeomorphic.
\end{corollary}
\begin{proof}
To prove that $\hat{\sg}_2(\hat{\Mg}_n)$ is open, take into account that
$\sg_2(\Mg_n)$ is open in  $\Delta^{3n+4}-\widetilde{H}_0$ and that
$\mbox{proj}_1:\Delta^{3n+4}-\widetilde{H}_0 \to \Cg_n$ is an open submersion.
Induce in $\hat{\Mg}_n$ the unique analytic structure making
$\hat{\sg}_2$ an analytic embedding.
Since $\hat{\sg}_2 \circ \mbox{proj}_2=\mbox{proj}_1 \circ \sg_2$ and
$\mbox{proj}_1:\sg_2(\Mg_n) \to \hat{\sg}_2(\hat{\Mg}_n)$ is an analytic
submersion, we infer that $\mbox{proj}_2$ is an analytic submersion.

For the second part of the corollary, let $\hat{\Sigma}_1$ and
$\hat{\Sigma}_2$ be two connected  components in $\hat{\Mg}_n,$ and suppose
that $\hat{\sg}_1(\hat{\Sigma}_1) \cap \hat{\sg}_1(\hat{\Sigma}_2) \neq
\emptyset.$ Take $(G,\mg_1) \in \mbox{proj}_2^{-1}(\hat{\Sigma}_1),$
$(G,\mg_2)\in \mbox{proj}_2^{-1}(\hat{\Sigma}_2).$ Denote by $\Sigma_1$
the connected component of $\Mg_n$ containing $(G,\mg_1),$ and analogously
define $\Sigma_2.$
Since $\sg_1(G,\mg_1)=\sg_1(G,\mg_2),$  Theorem \ref{th:structure} shows
the existence of $\tau \in {\cal P}_{n+1}$ such that
$\tau(\Sigma_1)=\Sigma_2,$ and so $\sg_1(\Sigma_1)=\sg_1(\Sigma_2).$
Moreover,  taking into account that $\mbox{proj}_3 \circ \sg_1=\hat{\sg}_1
\circ \mbox{proj}_2$  and $\mbox{proj}_2(\Sigma_j)=\hat{\Sigma}_j,$
$j=1,2,$  we get
$\hat{\sg}_1(\hat{\Sigma}_1)=\hat{\sg}_2(\hat{\Sigma}_2).$ In addition,
$\tau:\Sigma_1 \to \Sigma_2$ induces in a natural way an analytic
diffeomorphism $\hat{\tau}:\hat{\Sigma_1}\to \hat{\Sigma_2},$   which
concludes the proof.
\end{proof}
Following  Theorem \ref{th:structure}, it is natural to ask whether $\hat{\sg}_1:\hat{\Mg}_n \to \hat{\Gg}_n$ is
an analytic covering.
However the class of a marked graph  admitting symmetries has non-trivial isotropy group for the natural action $\hat{\mu}:{\cal
P}_{n+1} \times \hat{\Mg}_n \to \hat{\Mg}_n.$ Anyway, we can endow
$\hat{\Gg}_n-\mbox{proj}_3(\mbox{Sym}(\Gg_n))$ (where $\mbox{Sym}(\Gg_n)$
consists of the family of graphs with non-trivial symmetry group) with the
analytic structure making $\hat{\sg}_1:
\hat{\Mg}_n-\hat{\sg}_1^{-1}\big(\mbox{proj}_3(\mbox{Sym}(\Gg_n)) \big) \to
\hat{\Gg}_n-\mbox{proj}_3(\mbox{Sym}(\Gg_n))$ an analytic covering of
$(n+1)!$ sheets.

{\bf ISABEL FERNANDEZ, FRANCISCO J. LOPEZ, \newline
Departamento de Geometr\'{\i}a y Topolog\'{\i}a \newline
Facultad de Ciencias, Universidad de Granada \newline
18071 - GRANADA (SPAIN) \newline
e-mail:(first author) isafer@ugr.es, (second author) fjlopez@ugr.es
\vspace{0.2cm}
RABAH SOUAM, \newline
Institut de Mathématiques de Jussieu-CNRS UMR 7586\newline
Université Paris 7\newline
Case 7012\newline
2,place Jussieu\newline
75251 Paris Cedex 05, France
\newline
e-mail: souam@math.jussieu.fr}



\end{document}